\input amstex
\input amsppt1
\loadmsbm
\UseAMSsymbols
\NoBlackBoxes
\parindent=0pt
\magnification=1200

\font\pio=cmr10 scaled 1280

\def\empty{{\text{\rm \pio \o}}}

\def\cd{\Cal D}
\def\td{\widetilde{\Cal D}}

\def\cu{{\Cal U}}
\def\cv{{\Cal V}}
\def\cw{\Cal W}

\def\enpr{\quad \vrule height .9ex width .8ex depth -.1ex}

\def\bar{\overline}

\def\bre{\Bbb R}
\def\bna{\Bbb N}
\def\bin{\Bbb Z}

\def\defin#1.#2{\noindent \bigskip {${\bf Definition \ #1.#2}$}   }
\def\osser#1.#2{\noindent \bigskip {${\bf Remark     \ #1.#2}$}   }
\def\propo#1.#2{\noindent \bigskip {${\bf Proposition\ #1.#2}$}   }
\def\teore#1.#2{\noindent \bigskip {${\bf Theorem    \ #1.#2}$}   }
 \def\corol#1.#2{\noindent \bigskip {${\bf Corollary  \ #1.#2}$}   }
\def\lemma#1.#2{\noindent \bigskip {${\bf Lemma      \ #1.#2}$}   }
\def\vn#1{V^{(#1)}}

\def\llim{\mathop{\longrightarrow}}
\def\vo{V^{(0)}}
\def\von{V^{(1)}}
\def\rvo{\bre^{V^{(0)}}}

\def\vin#1{V^{(#1)}}
\def\rvin#1{\bre^{V^{(#1)}}}

\def\Lb{\Lambda}
\def\cl{{\Cal L}}
\def\ce{{\Cal E}}
\def\cf{{\Cal F}}
\def\cg{{\Cal G}}
\def\whcg{\widehat\cg}

\def\ta{\tilde A}
\def\qua{\quad $\,$}
\def\osc {\text{\rm Osc}}
\def\smad{\smallskip Proof.\ }
\def\witi{\widetilde}
\def\disp{\displaystyle}
\def\s,{\quad $\,$}
\def\veps{\varepsilon}
\def\bu{\Cal T}
\def\wt{\widetilde \cv}\font\pio=cmr10 scaled 1230
\def\empty{{\text{\rm \pio \o}}}
\def\rfbbkt{[1]}
\def\rfbft{[2]}
\def\rfhmt{[2]}
\def\rfhh{[3]}
\def\rfki{[4]}
\def\rfkib{[5]}
\def\rfl{[6]}
\def\rfm{[7]}
\def\rfmb{[8]}
\def\rfmc{[9]}
\def\rfmd{[10]}
\def\rfp{[11]}
\def\rfpb{[12]}
\def\rfpc{[13]}
 \def\rfpd{[14]}
\def\rfs{[15]}
\def\rfsb{[15]}
\def\rfst{[16]}

\centerline{\bf  Uniqueness of Eigenforms on  Fractals-II}
\medskip\centerline{\bf Roberto Peirone}

\smallskip\centerline{\it Universit\`a di Roma-Tor Vergata,Dipartimento di Matematica}          
\centerline{\it Via della Ricerca Scientifica, 00133 Roma, Italy}
\smallskip
\centerline{E.mail: peirone\@mat.uniroma2.it}

\bigskip{\bf Abstract.} 
{\it I give an explicitly verifiable
necessary and sufficient condition for the uniqueness of the 
eigenform on finitely ramified fractals, once an eigenform is known.
This improves the results of my previous paper \rfpd, where I gave
some necessary and some sufficient conditions, and with a relatively mild additional 
requirement on the known eigenform. The result of this paper can be interpreted as
a uniqueness result for self-similar energies on the fractal.}
 \smallskip
MSC: 31C25, 28A80,  Key Words: Dirichlet Forms, Fractals, Eigenforms
\bigskip

\centerline{\bf 1. Introduction}
\medskip 
This paper deals with analysis on  finitely ramified self-similar fractals.
 In order to define a
 self-similar fractal, we start with finitely many contractive 
 (i.e., having factor $<1$) similarities 
 $\psi_1,...,\psi_k$ in $\bre^{\nu}$ (more generally, in a complete metric space). 
 Then, we can say that there exists 
 a unique nonempty compact set $K$ in ${\bre}^{\nu}$ such that

$$K=\bigcup\limits_{i=1}^k \psi_i(K), \eqno (1.1)$$

 which we will call fractal or self-similar set generated by this set of contractions.
A tipical example of a self-similar fractal is the (Sierpinski) 
 Gasket where we have three  maps $\psi_i$, $i=1,2,3$, which are the rotation-free contractions 
 with factor ${1\over 2}$  that
  have  as fixed points the three vertices $P_i$ of a fixed equilateral triangle $T$,
   in formula 
  $\psi_i(x)=P_i+{1\over 2}(x-P_i)$.  Other examples are the Vicsek set,
  and the (Lindstr\o m) Snowflake.  The term finitely ramified means more or less
that the "copies" $\psi_i(K)$ of the fractal intersect only at finitely many points. An example
of an infinitely ramified fractal is the Sierpinski Carpet.
In the present paper, I will consider a subclass including Gasket, Vicsek set and Snowflake, of the set of the
 P.C.F. self-similar sets,  a class of 
 finitely ramified fractals introduced by J. Kigami  in \rfki. In fact I require a mild
additional hypothesis on the P.C.F.self-similar sets, as in other 
papers of mine, but also in other works (for example \rfst).
  A description of the general theory of P.C.F. self-similar sets with many examples can be 
    found in \rfkib, or also in \rfst.

\qua The problems discussed in this paper concern the construction of "energies" or in other words,
 Dirichlet forms, on fractals. This is a problem widely investigated in this area
and can be also interpreted as the construction of diffusions, or of harmonic structures
on the fractal. The self-similar Dirichlet forms, roughly speaking,
 are the "energies" $\ce$ for functions defined on $K$ with the property that
 
 $$ \ce(v)=\sum\limits_{i=1}^k r_i \ce(v\circ\psi_i)$$
 
 where the {\it weights} $r_i$ are suitable positive numbers.
 On the P.C.F. self-similar sets such  self-similar Dirichlet forms
are in correspondence
with the $r$-eigenforms, i.e., 
 the eigenfunctions of a
 special nonlinear operator $\Lb_r$ defined on the set of the  Dirichlet forms 
 on the {\it finite subset} of the
 fractal $\vo$, which is a sort of boundary of the fractal, for example the three
vertices of the triangle in the case of the Gasket. This correspondance is typical
of finitely ramified fractals, and this makes the problems concerning self-similar energies
on infinitely ramified fractals much harder.
 Such an operator $\Lb_r$ depends
 on the  set $r$ of weights $(r_i)_{i=1}^k$ placed on the "copies" $\psi_i(K)$ and  is called
the  {\it renormalization operator}.  

\qua  Natural problems in this context are the existence
and the uniqueness (up to a multiplicative constant) of the eigenform.
 Such  problems can
  be interpreted as the existence and the uniqueness of  a diffusion.
 Results concerning the existence of an eigenform with fixed weights on finitely ramified fractals
  has been proved given  in \rfl\ ,
   \rfs\ and in \rfmd.  In more recent papers, \rfhmt, \rfpb, \rfpc,
  conditions are given in order that there exists an $r$-eigenform corresponding
  to a suitable set of weights $r$. In particular, in \rfpc\ it is proved that
  this occurs on fractals with connected interior and on every fractal, 
 but at a suitable level, that is on the structure $\cf_n$ described in Section 4.
  
  \qua In this paper, I am interested 
  in the uniqueness of the $r$-eigenform (up
   to a multiplicative constant) on finitely ramified fractals. 
   On fractals which are not finitely ramified  there is a recent
   uniqueness result on Carpet-like fractals  in \rfbbkt. 
   On many usual (P.C.F.) fractals,  we have in fact 
   the uniqueness of the
   $r$-eigenform. The first  example of nonuniqueness was given by
  V. Metz  in \rfm\,where he proved 
   that on the Vicsek set we have infinitely many 
  normalized eigenforms. 
  The uniqueness of the
  eigenform has been investigated for example in \rfmb , \rfp , \rfs , 
  where  the uniqueness
  is proved on many specific fractals, but the 
  results given in such papers appear to be hardly applicable in  general situations. 
  Usually the methods there apply when we have three vertices or when  we
 require symmetry properties of the $r$-eigenform. 
 
 \qua A more general
 uniqueness criteria is given in \rfmc, where     Theorem 15
 provides a necessary and sufficient condition for the uniqueness 
 of the $r$-eigenform,
  once we know an $r$-eigenform $E$, in terms of the derivative 
 of the renormalization operator at $E$. Such
 a condition provides in many cases   an explicit, although
 usually rather technical, method to verify
  whether a given $r$-eigenform is unique. However,
   such a condition  does not appear to easily suggest a way to construct
 nontrivial cases of nonuniqueness.

  \qua 
  In my previous  paper \rfpd,  a method is described that allows us
to   prove  some sufficient and some necessary conditions
 for the uniqueness and can be used in rather general situations,
  and in many cases
 allows us to deduce the uniqueness or nonuniqueness 
 merely from the geometric
 structure of the fractal.  In particular the
  nonuniqueness cases include fractals that are very far from being a tree.
Moreover,  an example was given where
we have uniqueness or nonuniqueness of the eigenform, depending on the weights.
The general necessary conditions and sufficient conditions were given
 under  a relatively mild condition
 on the fractal denoted by (A), which is satisfied for example
 by fractals with connected interior.
 

\qua  However no explicit conditions given in \rfpd\  is
both necessary and sufficient. In this paper, I improve the results of \rfpd\ 
and give an explicit (that is, effectively verifiable) necessary and sufficient condition
for the uniqueness of a given eigenform. At a first step, I assume
Condition (A), and a modification of a condition given in \rfpd\
turns out to be necessary and sufficient for the uniqueness of a given eigenform.
Later, I remove condition (A), and give an explicit necessary and sufficient condition
for the uniqueness of the eigenform in the general case.

\qua When (A) holds, thus the eigenforms are positive, the methods
relies on the Perron-Frobenius Theory for the positive linear operator
(denoted by $T_{j;E;r}$) that sends a function $u$ on $\vo$ to the harmonic 
 extension of $u$ restricted to the cells containing the vertex $P_j$.
 Such an operator $T_{j;E;r}$ has a unique normalized positive eigenvector
 $\bar u_j$. The eigenform $E$ is not unique if and only if
 the set of $\bar u_j$ has two disjoint nonempty subsets that have a 
 special stability property. If we remove (A), we have preliminarly to split
 the complement of $P_j$ into components, that
  we denote by $C_{j,1},...,C_{j,m_j}$, and we associate
 an eigenvector $\bar u_{j,s}$ to every component $C_{j,s}$, $\bar u_{j,s}$
 being the eigenvector of a suitable power of $T_{j;E;r}$.
 In such a setting  the nonuniqueness amounts
 to the existence of two nonempty disjoint subsets of the sets of all
 $\bar u_{j,s}$ that are stable in a  sense similar to the previous.
 The proof in this general  case 
 is not a simple variation of the case where (A) holds, but requires some preparatory work.
 In particular, it is proved that the set positive coefficients of a given
 eigenform (which with (A) is the set of all coefficients) only
 depends on the fractal, thus is independent of the eigenform (Theorem 5.3).

\bigskip 
 \centerline{\bf 2. Notation and Preliminary Results}
 
 \medskip 
In the present Section,  I fix the general setting and
give the preliminary results.
First of all, I recall some notion on graphs.
We say that a pair  $(V,\cw)$ is a graph (or that {\it $\cw$ is  a graph on $V$})
    if  $\cw$ is  a set of subsets of $V$ having two elements (the {\it edges} of the graph).
    We say that a finite sequence  $(i_0,...,i_n)$ in $V$ 
    is a path (or more precisely a $\cw$-path) connecting $i_0$ and $i_n$ if
    $\{i_{h-1},i_{h}\}\in\cw$ for every $h=1,...,n$.
 Given $V'\subseteq V$, we say that the path is {\it in $V'$} if $i_h\in V'$ for every $h=1,...,n-1$.
 
\qua   We say that $(V,\cw)$ is connected 
    if any two points in $V$   
    can be connected by a path. 
    More generally, we say two points are connected in a subset
    $V'$ of $V$ if they are connected by  a path in $V'$. We say that $V'$
     is connected if every two points in $V'$ can be connected 
    by a path in $V'$. When $V'$ is not connected it can be split into its components.
    More precisely, if $x\in V'$ we call {\it component of $x$ in $V'$} the set of $y\in V'$ such
    that $x$ and $y$ are connected  by a path in $V'$. It can be easily verified that
    the components are connected and mutually disjoint and that
    $V'$ is the union of the  components in $V'$.

\qua  I now recall the basic definitions on fractals, following the same approach 
 as in \rfpd.    (a similar approach is discussed in
     \rfhh\ and in \rfp). Namely, start with a finite set $\Psi=\big\{\psi_1,...,\psi_k\big\}$ 
     of one-to-one 
   maps defined on a finite set $\vo=\big\{P_1,...,P_N\big\}$ 
   (not necessarily a subset of 
   $\bre^{\nu}$),  with $2\le N\le k$, and put
   
    $$\von=\bigcup\limits_{\psi\in\Psi}\psi(V^{(0)})\, .\eqno (2.1)$$ 
   
   We call $1$-cells the sets
    $V_i:=\psi_i(\vo)$ for $i=1,...,k$, and put 
    
    $$\cv:=\{i=1,...,k\}, \quad\cu=\{1,...,N\}, 
    \quad J:=\big\{\{j_1,j_2\}: j_1,j_2\in\cu, j_1\ne j_2\big\} .$$

  On $\cv$ we consider the graph $\cg_{\cf}$ whose edges are $\{i_1,i_2\}$ such that
    $i_1,i_2\in\cv$, $i_1\ne i_2$, and $V_{i_1}\cap V_{i_2}\ne\empty$. 
   We  require that for each $j=1,...,N$  
   
   $$\psi_j(P_j)=P_j,\qquad  P_j\notin\psi_i(V^{(0)})\quad\forall\,i\ne j
   \ \forall\,j=1,...,N,\eqno (2.2)$$
   $$ \cg_{\cf} \text{\rm\ is\ connected.}\eqno (2.3)$$
   
    We say that  $(\vo,\von,\Psi)$ is a {\it  fractal triple} and call it $\cf$.
    

  \qua   I will denote by
$\cd(\vo)$ or simply $\cd$ the set
of the  Dirichlet forms on $\vo$, invariant with respect to
an additive constant,  i.e., the set of the functionals $E$ 
from $\rvo$ into $\bre$ of the form

 $$E(u)=\sum\limits_{\{j_1,j_2\}\in J}
  c_{\{j_1,j_2\}}(E)\big(u(P_{j_1})-u(P_{j_2})\big)^2 $$
 
where the {\it coefficients} $c_{\{j_1,j_2\}}(E)$ (or simply 
$c_{\{j_1,j_2\}}$)   of $E$ are required to
be nonnegative. I will denote by
$\td(V)$ or simply $\td$ the set of the irreducible Dirichlet forms, i.e., $E\in\td$ if $E\in
\cd$ and moreover  $E(u)=0$
 if and only if $u$ is constant. 
 This amounts to the fact that 
the graph $\cg_0(E)$ is connected, where $\cg_0(E)$ is the graph on 
$\vo$ formed by all $\{P_{j_1}, P_{j_2}\}$ such that  $c_{\{j_1,j_{2}\}}(E)>0$.
 Every $E\in\cd$ 
  is uniquely determined by its coefficients. 
  I will say that $E$ is positive if all its coefficients  are positive.

\qua   Recall that   for every  $r\in W:= ]0,+\infty[^{\cv}$  ($r_i:=r(i))$
 the {\it renormalization operator} is defined as follows:
  for every $E\in\td$ and every $u\in\rvo$,  
    
 $$\Lb_{r}(E)(u)=\inf\Big\{S_{1;r}(E)(v), v\in\cl(u)\Big\},$$
 $$S_{1;r}(E)(v):=\sum\limits_{i=1}^k r_i E(v\circ\psi_{i}),\quad 
\cl(u):= \big\{v\in \bre^{V^{(1)}}: v=u
  \ \text{on}\ V^{(0)}\big\}.$$
 
 It is well known that $\Lb_r(E)\in\td$ and that
 the infimum is attained at a unique function $v:=H_{1,E;r}(u)$. 
 When $r\in W$, 
 an element $E$ of $\td$ is said to be an {\it $r$-eigenform} with 
 {\it eigenvalue} $\rho>0$ 
 if $\Lb_{r}(E)=\rho  E$. As this amounts to
 $\Lb_{{r\over \rho}}(E)= E$,  we could also require $\rho=1$.
  It is well-known that
 if there exist two $r$-eigenforms on the same fractal, 
then they have the same eigenvalue.
 I will say that there is $r$-existence if there exists an $r$-eigenform,
 $r$-uniqueness if the $r$-eigenform is unique up to a multiplicative constant,
 and I will similarly use the expressions $r$-nonexistence and 
 $r$-nonuniqueness. It is well-known that,
{\it   if $E$ is an $r$-eigenform, then $r_i>\rho$ for every
 $i\in\cu$ (not necessarily for every $i\in\cv$). }

 \qua  We say that $\von$ is $A$-connected,
    or that (A) holds, if 
    for every $j,j_1,j_2\in\cu$, with $ j_1\ne j\ne j_2$,
     there exists a path connecting $j_1$ to $j_2$  in $\cv\setminus \{j\}$.
    We say that $\von$ is $O$-connected or with connected interior if
    
   \qua a) The $1$-cells $V_j$, $j\in\cu$ are mutually disjoint,
    
    \qua b) $\cv\setminus\cu$ is connected.
    
    It is easy to verify that $O$-connected implies $A$-connected.
      Condition (A) is satisfied by many of the fractals discussed in literature.
       However, the tree-like Gasket, that is a fractal like the Gasket
   but  where two of the $1$-cells are separated, does not satisfy it.
   The interest of such notions is that 
   {\it  if (A) holds, then  every $r$-eigenform is positive.}

\qua  We define the linear operator $T_{i;E;r}$ from $\rvo$ into itself 
 by $T_{i;E;r}(u)=H_{1,E;r}(u)\circ \psi_i$, for $i=1,...,k$ and more generally
 $T_{1_1,...,i_n;E;r}:=T_{i_1;E;r}\circ\cdots\circ T_{i_n;E;r}$.  We also 
 define the standard operator $L_E$ (here $L$ stands for Laplacian)
 by

 $$L_{E}(u)(P_j)=\sum\limits_{h\ne j}  c_{j,h}
 \big(u(P_h)-u(P_j)\big)$$
 
  and we say that  $P_j$ 
   is an $E$-harmonic point for $u$ if $L_{E}(u)(P_j)=0$ and that $P_j$
 is an $E$-nonharmonic point for $u$ if $L_{E}(u)(P_j)\ne 0$. 
  Now, I recall the main results in this context.

 \bigskip
 {\bf Lemma 2.1.} {\sl Let $E\in\td$ and $u\in\rvo$. If
  $u$ is nonconstant, then $L_E(u)$ is positive
 at some point and negative at some point. If, moreover,  $E$ is positive,
  then $L_E(u)$ is positive at the minimum points of $u$
 and is negative at the maximum points of $u$.}
\smad
As $\cg_0(E)$ is connected and $u$ is not constant,  
there exist $j_1, j_2\in\cu$ such that
$c_{\{j_1,j_2\}}(E)>0$ and $u(P_{j_1})=\min u$, $u(P_{j_2})>\min u$. By  definition,
we have  $L_E(u)(P_{j_1})\ge c_{j_1,j_2}\big( u(P_{j_2})- u(P_{j_1})\big)>0$.
Similarly,  if $c_{\{j_1,j_2\}}(E)>0$ and $u(P_{j_1})=\max u$, $u(P_{j_2})<\max u$,
then $L_E(u)(P_{j_1})<0$. If $E$ is positive, we can apply the previous
considerations to {\it any} minimum/maximum point of $u$, and the
Lemma  is proved.  \enpr

 %
 %
 %
 %
 %

  \bigskip
 \qua If  $E,E'\in\td$, let $\lambda_+(=\lambda_+(E,E'))=\max {E'(u)\over E(u)}$,
  $\lambda_-(=\lambda_-(E,E'))=\min {E'(u)\over   E(u)}$, 
  where the maximum and the minimum are taken over all nonconstant $u$. Also, let
 
 $$ A^{\pm}(=A^{\pm}(E,E'))=\{u\in\bre^{V^{(0)}}: E'(u)=\lambda_{\pm}(E,E') E(u)\},$$
 $$ \tilde A^{\pm}(=\tilde A^{\pm}(E,E'))=\{u\in A^{\pm}: u \text{\ nonconstant}\}.$$
 
\qua We say that $X$ is a c-linear subspace of $\rvo$ if $X$ is a linear subspace
strictly  containing the constants. We say that a set $X$ is $T_{E;r}$-invariant if
 $u\in X$ implies $T_{i;E;r}(u)\in X$ for every $i\in\cv$.  We say that
 two c-linear subspaces $X, X'$ of $\rvo$ are almost disjoint
 if $X\cap X'$ is the set of the constants.

 \bigskip
  {\bf Proposition 2.2.} {\sl
 
\qua  i) $A^{\pm}(E,E')$ are c-linear subspaces and $T_{E;r}$-invariant, for any
 $E,E'$ $r$-eigenforms.
 
\qua  ii) $E'$  is a multiple of $E$ $\iff$ $\ta^+\cap\ta^-\ne\empty$.
 }
 \smad
 i) is in \rfs, Lemma 5.13 
 and ii) is in  \rfp, and is however trivial. \enpr


  \bigskip
  {\bf Lemma 2.3.} {\sl Suppose $\rho=1$. If there exist an $r$-eigenform
  $\witi E$ and a nonempty subset $H$ of $\td$, closed in $\td$ and invariant under $\Lb_r$ 
  not containing
  any multiple of $\witi E$, then we have $r$-nonuniqueness.   }
  \smad
  See \rfpd, Lemma 2.7. This depends on Theorem 4.22 in \rfp.
\enpr

 \bigskip
 {\bf Lemma 2.4} {\sl  Suppose $\witi E$ is an $r$-eigenform (not necessarily positive)
 with eigenvalue $1$. Suppose $X,X'$
 are  almost disjoint $T_{\witi E;r}$-invariant c-linear subspaces of $\rvo$. 
  Then the set  
 
 $$K_{t,\witi E, X,X'}:=\big\{ E\in\td:  E\le t \witi E,  E\le \witi
  E \ \text{on}\ X, E=t \witi E  \ \text{on}\ X' \big\}$$
  
  is closed in $\td$ and $\Lb_r$ invariant for every $t>1$.  }
 \smad
 I outline the idea. For the details see \rfpd, proof of Theorem 2.8.
 The set $K_t(:=K_{t,\witi E, X,X'})$ is  clearly closed.
 To prove that it is $\Lb_r$ invariant, we have to prove that for every $E\in K_t$, then
  $\Lb_r(E)\in K_t$ as well, and the most delicate
   point is to prove that   $\Lb_r(E)(u)= t \witi E(u) \ \ \forall\, u\in X'.$
     In order to prove it, it suffices to prove
  that 
  
  $$H_{1,\witi E;r}(u)=H_{1,E;r}(u)\quad \forall\, u\in X' .  \eqno (2.4)$$
 
Now, given $u\in X'$, note that 
 $\Big(S_{1;r}(E)-t S_{1;r}(\witi E)\Big)(v)=\sum\limits_{i=1}^k
r_i (E-t\witi E)\big(v\circ\psi_i)$
 by definition of $K_t$ is always non positive and attains the value $0$ at 
$v=H_{1,\witi E;r}(u)$ (as $H_{1,\witi E;r}(u)\circ\psi_i=T_{i;\witi E;r}(u)\in X'$).
It follows that $H_{1,\witi E;r}(u)$ is a stationary point for 
 $S_{1;r}(E)- t S_{1;r}(\witi E)$
 on $\cl(u)$, but, in view of a known result, 
  it is also a stationary point for
 $S_{1;r}(\witi E)$ on $\cl(u)$. Consequently,
 $H_{1,\witi E;r}(u)$ is a stationary point for $S_{1;r}(E)$ on  $\cl(u)$, thus,
 by the same result as above, it amounts
 to $H_{1,E;r}(u)$ and (2.4) is proved.  \enpr

 \bigskip
 {\bf Theorem 2.5.} {\sl Suppose  there exists a positive $r$-eigenform $\witi E$
 (this is true for every eigenform if (A) holds). Then
 there is $r$-nonuniqueness if and only if there exist
 two almost disjoint $T_{\witi E;r}$-invariant c-linear subspaces of $\rvo$.} 
\smad 
See \rfpd, Theorem 2.8.  \enpr

 \bigskip
 Suppose  $E$ is a positive $r$-eigenform,
  so that, or every $j\in\cu$
 $T_{j;E;r}$ is a positive linear operator 
  on $\rvo_j=:\big\{u\in\rvo: u(P_j)=0\big\}$. As a consequence,
  by the Perron-Frobenius theory (see Lemma 5.3 in \rfp)
   there exists a (unique) positive eigenvector
  $\bar u_j(=\bar u_{j,E;r})$  of $T_{j;E;r}$ on $\rvo_j$ of norm $1$ 
  with positive eigenvalue $l_j(=l_{j,E;r})$, and for every
  $u\in \rvo_j$ there exists $\pi_j(u)=\pi_{j,E;r}(u)$ such that 
 $\disp{{T_{j;E;r}^n(u)\over l_j^n}\llim\limits_{n\to+\infty} \pi_j(u) \bar u_j}$.
 Note that in \rfpd, I used the notation $\bar v_j$ in place of $\bar u_j$, but I think that
 $\bar u_j$ is more coherent with the other notation in this context.

 \bigskip
 {\bf Lemma 2.6.} {\sl 
 
 \qua i)  If   $E$ is an $r$-eigenform,
  then 
 $L_{E}(u)(P_j)=\big({{r_j\over \rho}}\big)^n L_{E}\big(T_{j;E;r}^n(u)\big)(P_j)$
  for every $P_j\in\vo$, every $u\in\rvo$ and every positive integer $n$.
 
\qua ii) If  $E$ is a positive
 $r$-eigenform, then $l_j=\rho r_j^{-1}\in]0,1[$.
 
\qua  iii) If  $E$ is a positive
 $r$-eigenform,  $X$ is a $T_{E;r}$-invariant c-linear
  subspaces of $\rvo$, $u\in X$ and $P_j$ is an 
  $E$-nonharmonic point for $u$, then $\bar u_j\in X$.}
  \smad
 See \rfpd, Lemmas 3.2 and 3.3.  Note that i) is a standard result. \enpr

\bigskip\smallskip
\centerline{\bf 3.  Uniqueness and Nonuniquenes Results for (A) connected fractals.} 

\bigskip
In \rfpd\ there are some necessary and some sufficient conditions 
for the uniqueness of the eigenform. In this Section I will introduce
a modification of one of the conditions which will turn out
to be necessary and sufficient for the uniqueness of the eigenform.

\qua If $E$ is a positive $r$-eigenform,
I say that {\it a subset $B$ of $\cu$ is $(E,r)$-$1$stable if, whenever $j\in B$ ,
$j'\in \cu\setminus B$
and  $i\in\cv$, then  $P_{j'}$
is an $E$-harmonic point for $T_{i;E;r}(\bar u_j)$.}

\qua If $E$ is an $r$-eigenform, {\it I say that a subset $B$ of $\cu$ is $(E,r)$-hstable
 if the set 

\centerline{$\bu_{B,E}:=
\big\{u\in\rvo:$ every point of $\cu\setminus B$ is $E$-harmonic for $u\big\}$}

 is  $T_{E;r}$-invariant.}

\qua The results connecting such notions with the uniqueness of the eigenform
are the following 

\qua {\it  (Lemma 4.4 in \rfpd) If $E$ is a positive $r$-eigenform, and every two   
 nonempty $(E,r)$-$1$stable sets  
 are not disjoint, then there is $r$-uniqueness.}

 \qua {\it (Lemma 4.5 in \rfpd) If $E$ is a positive $r$-eigenform and there exist two  disjoint
 $(E,r)$-hstable subsets $B$ and $B'$ of $\vo$  having at least two elements,
 then there is $r$-nonuniqueness.}

\qua In \rfpd, in fact, the term $(E,r)$-stable instead of $(E,r)$-$1$stable was used.
I here prefer to reserve the name $(E,r)$-stable to the new notion
I am introducing here, which will provide a necessary and sufficient condition
for the uniqueness and to see the previous notion as the $1$-case of it.
Namely, if $E$ is a positive $r$-eigenform,

\qua {\it I say that a subset $B$ of $\cu$ is $(E,r)$-stable if, for every $j\in B$,
$j'\in \cu\setminus B$, every $n\in\bna$
and every $i_1,...,i_n\in\cv$,  then  $P_{j'}$ is an 
 $E$-harmonic point for $T_{i_1,...,i_n;E;r}(\bar u_j)$.}
 
 \qua
 At first glance, one could think that, in order to check whether a set is $(E,r)$-stable,
 we have to verify infinitely many conditions, but, in fact,
 we can reduce the problem to only finitely many conditions. In fact,
 given $B\subseteq \cu$, let $X_{B, n}$ be the linear space generated by

$$\big\{T_{i_1,...,i_h;E;r}(\bar u_j):i_1,...,i_h=1,...,k, h\le n, j\in B\big\}.$$

Then we clearly have
 $X_{B,n}\subseteq X_{B,n+1}$. Moreover, if the equality holds, then
 $X_{B,m}=X_{B,n}$ fo every $m\ge n$. In fact, 
 $T_{i_1,...,i_{n+2};E;r}(\bar u_j)=T_{i_1;E;r}\big(T_{i_2,...,i_{n+2};E;r}(\bar u_j)  \big)$ 
 and as, by hypothesis $T_{i_2,...,i_{n+2};E;r}(\bar u_j) $ is a
  linear combination of elements in $X_{B,n}$, then
  $T_{i_1,...,i_{n+2};E;r}(\bar u_j)$ is a linear combination of elements of
  $T_{i_1;E;r}\big(X_{n}\big)$  thus is in $X_{B,n+1}$.
  Now, as $\rvo$ has dimension equal to $N$ and $X_{B,0}$
  has positive dimension, then there exists
  $n\le N-1$ such that $X_{B,n}=X_{B,n+1}$, hence
  $\bigcup\limits_{n=0}^{+\infty}X_{B,n}=X_{B,N-1}$.
  
  \qua It follows that a subset $B$ of $\cu$ is $(E,r)$-stable if and only if, for every $j\in B$, 
$j'\in\cu\setminus B$,  $n\le N-1$, $i_1,...,i_n\in\cv$, then $P_{j'}$ is 
 $E$-harmonic  for $T_{i_1,...,i_n;E;r}(\bar u_j)$.
In fact, by the argument above, if a point is $E$-harmonic for
every function  of the form $T_{i_1,...,i_n;E;r}(\bar u_j)$, $n\le N-1$, then 
it is also $E$-harmonic for every function of the form 
$T_{i_1,...,i_n;E;r}(\bar u_j)$, $n\in \bna$.
Hence, we can verify with finitely many calculations what subsets of $\cu$ are 
$(E,r)$-stable.

 \bigskip
{\bf Theorem 3.1}. {\sl If $E$ is a positive $r$-eigenform, 
then we have $r$-uniqueness if and only if every two nonempty 
$(E,r)$-stable sets are not disjoint.}
\smad
 Suppose every two nonempty $(E,r)$-stable sets are not disjoint, and prove 
 that there is $r$-uniqueness. By Theorem 2.5, we have to prove that, given two 
 $T_{E;r}$-invariant
c-linear subspaces of $\rvo$ $X_1$ and  $X_2$ , then 
 $X_1$ and $X_2$ are not almost disjoint. Let
 $B_h:= \{j:\bar u_j\in X_h\}$  for $h=1,2$.
 
 \qua If $j \in B_h$, then $T_{i_1,...,i_n;E,r}(\bar u_j) \in X_h$ as $X_h$ is 
$T_{E;r}$-invariant, thus, by Lemma 2.6 iii)
$B_1$ and $B_2$ are nonempty and $(E,r)$-stable. 
Thus, there exists $j \in B_1 \cap B_2$, that is, 
$\bar u_j \in X_1 \cap X_2$, so that  $X_1$ and $X_2$ are not almost disjoint.

\qua For the converse, suppose there exist two nonempty 
disjoint $(E,r)$-stable  subsets of $U$ $B_1$ and $B_2$. 
Then for every $i_1,...,i_n \in\cv$ and $h \in B_i$, $j\in U
 \setminus B_i$, $P_j$ is $E$-harmonic for $T_{i_1,...,i_n;E;r}(\bar u_h)$. Now put
 
 $$X_i := \big\{u\in\rvo : P_j\,  \text{is}\,  E\text{-harmonic\  for}\,
 T_{i_1,...,i_n;E;r}(u) \, \forall\, j \in U \setminus  B_i,  \forall\, n\in\bna,\,
  i_1,...,i_n \in \cv\big\}.$$
 
Then, for every $i = 1,2$  $X_i$ contains $\bar u_h$ with 
 $h \in B_i$, so that $X_i$ is a c-linear subspace of $\rvo$. Next, $X_1$ and $X_2$
  are almost disjoint. In fact, if $u \in X_1 \cap X_2$, then every $j \in \cu$
either lies in $\cu \setminus B_1$ or in $\cu \setminus B_2$, as $B_1$ and $B_2$ are 
  supposed to be disjoint. Hence, by hypothesis, considering
 $n=0$ in the definition of $X_i$, every $P_j\in \vo$ is $E$-harmonic  for $u$, 
  hence, by Lemma 2.1,
   $u$ is constant. Finally, $X_1$ and $X_2$ are clearly $T_{E,r}$-invariant.
   By Theorem 2.5 we have $r$-nonuniqueness.  \enpr

 \bigskip
 {\bf Remark 3.2.}
 We will now see that Lemmas 4.4 and 4.5 in \rfpd\ can be seen as direct 
 consequences of 
 Theorem 3.1 here, hence all results in \rfpd\ can be deduced from it.
 Clearly, an $(E,r)$-stable set is $(E,r)$-$1$stable, and we are going to see
 that every $(E,r)$-hstable set $B$ having at least two elements is
 $(E,r)$-stable. Let $j\in B$ and take $j_1\in B$ with
 $j_1\ne j$. Let $u\in \rvo$ be a function
 that attains the value $1$ at $P_{j_1}$ and $0$ at  $P_j$ and
 is $E$-harmonic at the other points. Thus, $u\in \bu_{B,E}$. As 
 $P_j$ is $E$-nonharmonic for $u$ by Lemma 2.1,
 and $\bu_{B,E}$ is by definition, $T_{E,r}$-invariant,
 and is clearly a c-linear subspace of $\rvo$, by Lemma 2.6 iii) $\bar u_j\in \bu_{B,E}$.
 Moreover, by the $T_{E,r}$-invariant again of $\bu_{B,E}$, if 
 $j'\in\cu\setminus B$, then $P_{j'}$ is $E$-harmonic for
 $T_{i_1,...,i_n;E;r}(\bar u_j)$ for every $i_1,...i_n\in\cv$. In conclusion,
 $B$ is $(E,r)$-stable, as required.  \enpr
 \bigskip

 \centerline{\bf 4. Nonpositive eigenforms.}

\bigskip

\qua Given a fractal triple $\cf:=(\vo,\von,\Psi)$, $\Psi=\{\psi_1,...,\psi_k\}$, we 
can define a related fractal set $K_{\cf}$ and a set
$\vin {\infty}$ that satisfy $\vo\subseteq\vin {\infty}\subseteq K_{\cf}$. Moreover,
 for any $i=1,...,k$ we can define
a one-to-one map from $\vin{\infty}$ into itself extending $\psi_i$, which we will still call
$\psi_i$. As we will never use  the set $K_{\cf}$, I will not describe it. 
On the contrary, I will describe $\vin {\infty}$ and the derived $n$-triples as they 
are essential for the following. However, I will merely recall the properties
and will not give the details, which are standard (see, e.g. \rfpc, Section 5). We put

$$\psi_{i_1,...,i_n}:=\psi_{i_1}\circ\cdots\circ\psi_{i_n},
\quad A_{i_1,...,i_n}=\psi_{i_1,...,i_n}(A)\quad \forall\,  A\subseteq K\, ,
$$
$$\vn{n}:=\bigcup\limits_{i_1,i_2,...,i_n=1}^k V_{i_1,...,i_n}, 
\quad \vn{\infty}=\bigcup\limits_{n=1}^{\infty}\vn{n}\, .$$

 The sets $V_{i_1,...,i_n}:=\vo_{i_1,...,i_n}$ are called {\it $n$-cells}. 
 Moreover, the sets are constructed in such  a
way that, if $(i_1,...,i_n)\ne(i'_1,...,i'_n)$, then 

$$\vin{\infty}_{i_1,...,i_n}\cap \vin{\infty}_{i'_1,...,i'_n} = V_{i_1,...,i_n}\cap
V_{i'_1,...,i'_n}\, . \eqno{(nesting\ axiom)}$$

As a consequence,  if  $(i_1,...,i_n)\ne(i'_1,...,i'_n)$ and 
$\psi_{(i_1,...,i_n)}(Q)=\psi_{(i'_1,...,i'_n)}(Q')$ with $Q,Q'\in K$, then $Q,Q'\in\vo$.
 Given the fractal triple $\cf$, we also define a related $n$-fractal triple
$\cf^n$ by

$$\cf^n:=(\vo,\vn{n},\Psi^n),\quad \Psi^n:=\big\{\psi_{i_1,...,i_n}: i_1,...,i_n\in\cv \big\}.$$

Here, the map $\psi\in\Psi_n$ satisfying
$\psi(P_j)=P_j$ is $\psi_{j^{(n)}}$ for every $j\in\cu$. 
Here, for every $i\in\cv$ $i^{(n)}$ denotes $i$ repeated $n$ times. 
Thus, we can define $\cf_n$-$\cv=\cv_n$, $\cf_n$-$\cv=\{ j^{(n)}:j\in\cu\}$
and $P_j^{(n)}=P_j$.
Moreover, as $\cf_n$ is a fractal triple, and consequently 
it satisfies (2.2), if $i_1,...,i_n\in\cv$ and $j,h\in\cu$,
then we have $\psi_{i_1,...,i_n}(P_h)=P_j$ if and only if $i_1=\cdots=i_n=j=h$.
Let 

$$W_n=\{r:\cv^n\to ]0,+\infty[\},$$ 

so that $r\in W_n$ can be written as
$r=\big(r_{i_1,...,i_n}\big):i_1,...,i_n\in\cv$. In particular, if $r\in W$ we define
$r^n\in W_n$ by $r^n_{i_1,...,i_n}=r_{i_1}\cdots r_{i_n}$.
 If $E\in\td$ and $r\in W_n$
we denote by $\cf_n$-$\Lb_r(E)$
the form $\Lb_r(E)$ in the triple $\cf_n$. This can be seen in the following way

$$\cf_n\text{-}\Lb_{r}(E)(u)=\inf\Big\{S_{n;r}(E)(v), v\in\cl(n,u)\Big\}$$
 $$S_{n;r}(E)(v):=\sum\limits_{i_1,...,i_n=1}^k r_{i_1,...,i_n} 
 E(v\circ\psi_{i_1,...,i_n}),\quad 
\cl(n,u):= \big\{v\in \rvin n: v=u
  \ \text{on}\ V^{(0)}\big\}.$$
 
The infimum is attained at a unique function $v:=H_{n,E;r}(u)$.  When $r\in W$, we write
$\Lb_{r^n}(E)$ short for $\cf_n\text{-}\Lb_{r^n}(E)$ as we know that $r^n$ lies in $W_n$.
We will occasionally  use the following generalization. Suppose $V'$ is a nonempty subset of  $\vin n$,
$r\in W_n$ and $w:V'\to\bre$. Then, we will denote by $H_{\vin n,V'}\big(S_{n,r}(E)\big)(w)$
the unique function from $\vin n$ into $\bre$ that amounts to $w$ on $V'$ and minimizes
$S_{n,r}(E)$. The following is well-known.

 \bigskip
 {\bf Lemma 4.1} {\sl 
 
 i) $\Lb_{r^n}(E)=\Lb_r^n(E)$ for every $E\in\td$ and $r\in W_1$.
 
 ii) $\cf_n$-$T_{i_1,...,i_n;E;r^n}(u)=T_{i_n;E;r}\circ\cdots\circ T_{i_1;E;r}(u)$
for every $u\in\rvo$,  if $E\in\td$ is an $r$-eigenform.} 
  \enpr

\bigskip
\qua  Note that in view of Lemma 4.1 i), if $E$ in an $r$-eigenform in
$\cf$, then $E$ is an $r^n$-eigenform in $\cf_n$, and such a fact will be used
implicitly in Section 6, as we will always assume there that
$E$ is an $r$-eigenform, and in some points we will have to consider $E$
in $\cf_n$. We are now going to define some graphs 
useful for the rest of the paper.
First, we define graphs only depending on $\cf$.
Let $\cg_1$ be the graph on $\von$ whose edges are the sets of the form

$$\big\{\psi_i(P_{j_1}), \psi_i(P_{j_2})\big\},\ \  i\in\cv, \{j_1,j_2\}\in J,$$

 and
$\cg'_1$ be the graph on $\von$ whose edges are the sets of the form

$$\big\{\psi_i(P_{j_1}), \psi_i(P_{j_2})\big\},\ \  i\in\cv\setminus\cu, \{j_1,j_2\}\in J.$$
 Let $\witi \cg$ be the graph on $\vo$ defined by the set of $\{P_{j_1}, P_{j_2}\}$ such that
$\{j_1,j_2\}\in J$, and there exist $Q_1\in V_{j_1}$ and $Q_2\in V_{j_2}$ 
 that are connected in the graph $\cg'_1$.
Next, given $E\in\td$ we define some graphs related to  $E$.  
Let $\cg_n(E)$ be the graph on $\vin n$ defined by
 
 $$\big\{\{\psi_{i_1,...,i_n}(P_{j_1}), \psi_{i_1,...,i_n}(P_{j_2})\}: \{j_1,j_2\}\in J, 
c_{j_1,j_2}(E)>0\big\}.$$

\qua  Given a graph $\cg$ on $\vo$, let $S(\cg)$ be the graph
 on $\von$ defined by the set of 

$\big\{\psi_i(P_{j_1}),\psi_i(P_{j_2})\big\}$ such that
 $\{P_{j_1},P_{j_2}\}\in \cg$.
Note that $S(\cg_0(E))=\cg_1(E)$. Let
  $\Lb(\cg)$ be the graph on $\vo$ defined by the set of 
 $\{P_{j_1}, P_{j_2}\}$ such that  $\{j_1,j_2\}\in J$, and 
$P_{j_1}$ and $P_{j_2}$  are connected in $\von\setminus\vo$
 by a path in $S(\cg)$.  We will write $\cf_n$-$\Lb(\cg)$, $\cf_n$-$\witi\cg$ 
 and so on to denote that 
the graphs are defined in the triple $\cf_n$.
 The statements of the following lemma are either well known or trivial.

 \bigskip
 {\bf Lemma 4.2} {\sl

i) If $\cg$ and $\cg'$ are graphs on $\vo$ 
and $\cg\subseteq\cg'$, then $\Lb(\cg)\subseteq\Lb(\cg')$. 

 ii) The graphs $\cg_n(E)$, $E\in\td$, and $\witi\cg$ are connected. Moreover,
if $\cg$ is a connected graph on $\vo$, so is also  $\Lb(\cg)$.
 
 iii) $\Lb\big(\cg_0(E)\big)=\cg_0\big(\Lb_r(E)\big)$ for every $E\in\td$ and $r\in W$.
 
 iv) If $E$ is an $r$-eigenform, then  $\Lb\big(\cg_0(E)\big)=\cg_0(E)$.
 
 v) $\Lb(\cg)\supseteq\witi\cg$ for every connected graph $\cg$ on $\vo$. } \enpr

 \bigskip
 It follows that $\Lb^n(\witi\cg)$ is increasing in $n$, hence there exists
 $n_1$ such that $\Lb^n(\witi\cg)=\Lb^{n_1}(\witi\cg)$ for every $n\ge n_1$.
 Put $\whcg=\Lb^{n_1}(\witi\cg)=\bigcup\limits_{n=0}^{+\infty}
 \Lb^{n}(\witi\cg)$. Note that clearly $\Lb(\whcg)=\whcg$.

 \bigskip
 {\bf Corollary 4.3} {\sl If $E\in\td$ is an $r$-eigenform, then 
 $\cg_0(E)\supseteq\whcg$.}
 \smad
 We have $\cg_0(E)\supseteq\witi\cg$ by Lemma 4.2 v), and by recursion,
 in view of Lemma 4.2, i) and iv), $\cg_0(E)\supseteq\Lb^n(\witi\cg)$ 
 for every $n$, in particular for $n=n_1$. \enpr

\bigskip
The following Lemma provides some standard variants 
of the well-known {\it maximum principle},
which corresponds to Lemma 4.4 i) in the case $C=\cv$.

 \bigskip
{\bf Lemma 4.4.} {\sl Let $E\in\td$.

\qua i) Suppose $V'$ is a nonempty subset of $\von$,
$w:V'\to\bre$, $C$ is  a nonempty connected
subset of $\cv$,  let $v:=H_{\von,V'}(S_{1;r}(E))(w)$, and put 
$\partial=V'\cup \big(V(C)\cap V(\cv\setminus C)\big)$. Then
  $v$ attains its
extrema (maximum and minimum) 
on $V(C)$ at points  in $\partial\cap V(C)$.
 
\qua ii) Given $Q\in\von\setminus\vo$, $u\in\rvo$, let $v=H_{1;E;r}(u)$ and let

 $$R_E(Q):=\big\{ P\in\vo: P\text{\ and\ }Q\text {\ are\ connected\ in }\ 
 \big(\von\setminus\vo, \cg_1(E)\big)\big\}$$
 
 Then $\min\limits_{R_E(Q)} u\le v(Q)\le \max\limits_{R_E(Q)} u$
 and the inequalities are strict unless $u$ is constant on $R_E(Q)$.
 
 \qua iii) If $u\ge 0$, then $T_{j;E;r}(u)\ge 0$.
 }
 \smad (Hint). i) Suppose that for some $Q\in V(C)\setminus\partial$, $v(Q):=M$  is the
maximum (resp. minimum) of $v$ on $V(C)$.
Then, $v(Q')=v(Q)$ if $\{Q',Q\}\in \cg_1(E)$, and such a point $Q'$ belongs to $V(C)$ by 
the definition of $\partial$. Hence, either i) for maximum (resp. minimum)
 holds or $v(Q')=M$ for every $Q'\in V(C)$.
But, if $C=\cv$, this implies $v(Q')=M$ for every $Q'\in V'$,  if
$C\subsetneq \cv$, this implies $v(Q')=M$ for some $Q'\in V(C)\cap V(\cv\setminus C)$, and
i)  for maximum (resp. minimum) holds again. ii) is well-known (e.g., it is in  \rfp,
Prop. 4.10) and can be proved like the maximum principle,
and iii) is an immediate consequence of i) with $C=\cv$ and
$V'=\vo$. \enpr

 \bigskip

 \centerline{\bf 5. Graph of an $r$-eigenform.}

 \bigskip
 {\bf Lemma 5.1.} {\sl Let $w:V'\to\bre$, $V':=\vo\setminus  \{P_{j_1}\}$ be defined by
 $w=\chi_{\{P_{j_2}\}}$. Let $v=H_{\von,V'}(S_{1;r}(E))(w)$.
 Then $v(P_{j_1})=\disp{ {c_{j_1,j_2}(\Lb_r(E))\over\sum\limits_{j\ne j_1}
 c_{j_1,j}\big(\Lb_r(E)\big) }}$. }
 \smad
 As $v|_{\vo}$ minimizes $\Lb_r(E)$ among the functions $v':\vo\to\bre$ such that
 $v'=w$ on $V'$, then the value $x:=v(P_{j_1})$  minimizes
 
 $$\sum\limits_{j\ne j_1} c_{j,j_1} \big(\Lb_r(E)\big)  \big(x-v(P_j)\big)^2
 =c_{j_1,j_2}(\Lb_r(E))(x-1)^2
 +\sum\limits_{j_1\ne j\ne j_2} c_{j,j_1} \big(\Lb_r(E)\big) x^2.$$

  Thus taking the derivatives we have
  
  $$x \sum\limits_{j\ne j_1} c_{j,j_1} \big(\Lb_r(E)\big)=c_{j_1,j_2}(\Lb_r(E)). \enpr $$

\bigskip
{\bf Lemma 5.2.} {\sl If $\{j_1, j_2\}\in J$ and $j_1^{(n)}$ and $j_2^{(n)}$ are connected in
$\vin n\setminus \big\{ j^{(n)}: j\in\cu\setminus \{j_1,j_2\}\big\}$ 
in the graph $\cg_{\cf_n}$, then
$\{P_{j_1},P_{j_2}\}\in \whcg$.}
\smad
Let $E\in\td$ be so that $\cg_0(\witi E)= \witi\cg$. Then

$$\{P_{j_1},P_{j_2}\}\in \cf_n\text{-}\witi G\subseteq
\cf_n\text{-}\Lb (\cg_0)= \cg\big(\Lb_{r^n}(\witi E)\big)=
 \cg\big(\Lb_r^n(\witi E)\big)=\Lb^n(\witi\cg)\subseteq \whcg. \enpr$$

 \bigskip
 {\bf Theorem 5.3.} {\sl If $E\in\td$ is an $r$-eigenform, then 
 $\cg_0(E)=\whcg$.}
 \smad
 In view of Corollary 4.3 it suffices to prove the inclusion "$\subseteq$".
 Thus, it suffices to prove that given $\{P_{j_1},P_{j_2}\}\notin \whcg$, then
 $c_{j_1,j_2}(E)=0$.  Let  $w_n$ and $v_n$ be as in Lemma 5.1 where we take $\cf_n$
in place of $\cf$ and $r^n$ in place of $r$.
 We are going to prove that 
 
 $$v_n(P_{j_1})\llim\limits_{n\to +\infty}0.\eqno (5.1)$$
 
 As  $E$ is an $r$-eigenform, in view of Lemma 5.1 we have
 
 $$v_n(P_{j_1})=\disp{ {c_{j_1,j_2}(\Lb_{r^n}(E))\over\sum\limits_{j\ne j_1}
 c_{j_1,j}\big(\Lb_{r^n}(E)\big) }}=
\disp{ {c_{j_1,j_2}(\Lb_r^n(E))\over\sum\limits_{j\ne j_1}
 c_{j_1,j}\big(\Lb_r^n(E)\big) }} =\disp{ {c_{j_1,j_2}(E)\over\sum\limits_{j\ne j_1}
 c_{j_1,j}(E)  }}$$
 
 and  by (5.1) we deduce $c_{j_1,j_2}(E)=0$, as required.
 By the maximum principle we have $0\le v_n\le 1$ on $\von$. Let 
 now $j\in\cu\setminus\{j_1,j_2\}$.
 By Lemma 4.1 ii) we have 
 
$$v_n\circ\psi_{j^{(n)}}= T_{j;E;r}^n(v_n|_{\vo}).$$ 

On the other hand, by the maximum
 principle and a compactness argument there exists $\bar\alpha\in ]0,1[$
 such that $\osc(T_{j;E;r}(u))\le\bar \alpha \osc (u)$. Hence,
 for every $\veps>0$ there exists $\bar n$ such that for every such $j$ and $n\ge \bar n$
we have $\osc (v_n\circ\psi_{j^{(n)}})\le\veps$. 
As $v_n(P_j)=0$ by definition, we then have

$$ v_n\circ\psi_{j^{(n)}}\le \veps. \eqno (5.2)$$

 Let  $C$ be the component of $j_1^{(n)}$ in
 $\{(i_1,...,i_n)\}\setminus\{j^{(n)}:j_1\ne j\ne j_2\}$
and take $n\ge \bar n$. As $\{P_{j_1},P_{j_2}\}\notin \whcg$,
by Lemma 5.2 $j_2^{(n)}\notin C$, hence $V'\cap \partial =\empty$.
  By Lemma 4.4 i), $v_n$ takes the maximum
  on $V(C)$ at points in $\partial$ that in our case is contained
  in $\bigcup\limits_{j_1\ne j\ne j_2} V_{j^{(n)}}$.
 Thus, by (5.2) we have $0\le v_n(Q)\le\veps$ for every $Q\in V(C)$, in particular
 for $Q=P_{j_1}$ and (5.1) is proved. \enpr

 \bigskip
 \centerline{\bf 6. Components of $\vo\setminus\{P_j\}$.}

 \bigskip
  Given $j\in\cu$, let $C_{j,1},...,C_{j,m_j}$ be the components of
 $\vo\setminus\{P_j\}$ in the graph $\whcg$, which,
 in view of  Theorem 5.3 amounts to $\cg_0(E)=\Lb\big(\cg_0(E)\big)$,
 whenever $E\in\td$ is an $r$-eigenform. 
 In this Section we will suppose that $E$ is an $r$-eigenform.
 We will not take care of determining what results are also valid
 for a generic $E\in\td$.  When $j'\in\cu\setminus\{ j\}$, let 
 
 $$L_j(P_{j'}):=\big\{ P_h\in\vo\setminus\{P_j\}:
P_{j'} \ \text{ is\ connected\ to \ }  \psi_j(P_h)  \text{\ in\ } \big(
\von\setminus\vo, S(\whcg) \big)\big\}.$$
 
  More generally,
 we put $L_j(B)=\bigcup\limits_{P_{j'}\in B} L_j(P_{j'})$ when 
 $B\subseteq \vo\setminus \{P_j\}$. In Lemma 6.1, we will prove some basic
properties of the sets  $L_j(B)$. In particular,  Lemma 6.1 i) provides
 an equivalent formulation of $L_j(P_{j'})$ in terms of an $r$-eigenform
  $E$. Such  a formulation will be used in the following without mention.

 \bigskip
 {\bf Lemma 6.1}
{\sl  Let $E$ be an $r$-eigenform. Then

\qua i) Given $P_h\in\vo$, then  $P_h\in L_j(P_{j'})$ if and only if
$T_{j;E;r}(\chi_{P_{j'}})(P_h)>0$.

\qua ii) Given $\{j,j'\}\in J$,
we have $L_j(P_{j'})\ne\empty$ if and only if $\{P_j,P_{j'}\}\in \whcg$.

\qua iii) If $j,j_1,j_2$ are mutually different elements of
$\cu$ and $L_j(P_{j_1})\cap L_j(P_{j_2})\ne\empty$ 
then $\{P_{j_1},P_{j_2}\}\in\whcg$, in particular 
$P_{j_1}$ and  $P_{j_2}$ are in the same $j$-component.

\qua iv) For every $\{j, j'\}\in\ J$,  $L_j(P_{j'})$ is either empty or a $j$-component.

\qua v) There exists a bijiection $\beta=\beta_j$ from $\{1,...,m_j\}$ into itself 
such that $L_j(C_{j,s})=C_{j,\beta(s)}$.}
\smad
i) We have $T_{j;E;r}(\chi_{P_{j'}})(P_h)=H_{1;E;r}(\chi_{P_{j'}})
\big(\psi_j(P_h)\big)$, thus,
by Lemma 4.4 ii),  $T_{j;E;r}(\chi_{P_j'})(P_h)>0$
if and only if $h\ne j$ and $P_{j'}$ is connected to 
$\psi_j(P_h)$  in $(\von\setminus\vo, \cg_1(E))$.  Now recall 
that $\cg_1(E)=S\big(\cg_0(E)\big)$ and $\cg_0(E)=\whcg$.

ii) Recall that $\Lb(\whcg)=\whcg$, hence, by the definition
of $\Lb(\whcg)$, $\{P_j, P_{j'}\}\in\whcg$ if and only if
 $P_j$ and $P_{j'}$ are connected in $\big(\von\setminus\vo,
S(\whcg)\big)$. This in turns, clearly amounts to 
the  existence of  $P_h$ with $h\ne j$ such that 
$P_{j'}$ is connected to  $\psi_j(P_h)$  in 
$\big(\von\setminus\vo, S(\whcg)\big)$,
and ii) follows at once.

iii)   Take  $P_h\in  L_j(P_{j_1})\cap L_j(P_{j_2})$.  Then $h\ne j$, 
so that $\psi_j(P_h)\in\von\setminus\vo$. Moreover
both $P_{j_1}$ and $P_{j_2}$ are connected
to $\psi_j(P_h)$,  hence $P_{j_1}$ and $P_{j_2}$ are connected, in
$\big(\von\setminus\vo, S(\whcg)\big)$. It follows that
$\{P_{j_1},P_{j_2}\}\in\Lb(\whcg)=\whcg.$

iv) and v) 
If $P_h$ and $P_{h'}$ lies in the same $j$-component, then 
$\psi_j(P_h)$ and $\psi_j(P_{h'})$ are connected
in $\big(\von\setminus\vo, S(\whcg)\big)$.
Hence $P_h\in L_j(P_{j'})$ if and only if
$P_{h'}\in L_j(P_{j'})$. Thus,  $L_j(P_{j'})$ is the union
of $j$-components, hence,  for every $s=1,...,m_j$, so is $L_j(C_{j,s})$.
On the other hand as $\whcg$ is connected, by the definition 
of $j$-component, there exists $P_{j'}\in C_{j,s}$ such that
$\{P_j,P_{j'}\}\in\whcg$, hence, by ii), $L_j(C_{j,s})\ne\empty$.
Therefore, there exists a nonempty subset $\Gamma_{j,s}$
of $\{1,...,m_j\}$ such that

$$L_j(C_{j,s})=\bigcup\limits_{s'\in \Gamma_{j,s}} C_{j,s'}.$$

On the other hand, the nonempty sets $\Gamma_{j,s}$,
$s=1,...,m_j$, are mutually disjoint by iii), hence
they are all singletons, and v) is proved.
Finally, given $P_{j'}\in\vo\setminus\{P_j\}$, let $C_{j,s}$ be the
$j$-component containing it. Then
$L_j(P_{j'})$ is the union of $j$-components, and is
contained in $L_j(C_{j,s})=C_{j,\beta(s)}$, thus  either is empty 
or amounts to $C_{j,\beta(s)}$, and iv) is proved. \enpr

 \bigskip
 When $u\in\rvo$ we have $u=\sum\limits_{P_{j'}\in\vo} u(P_{j'})\chi_{P_{j'}}$, hence
 for every $j\in\cv$, $r\in W$ and $n\in\bna$ we have
 
 $$T_{j;E;r}^n(u)=\sum\limits_{P_{j'}\in\vo}  u(P_{j'})  T_{j;E;r}^n (\chi_{P_{j'}})
\eqno (6.1) $$

  \bigskip
 {\bf Lemma 6.2} {\sl Suppose $E$ is an $r$-eigenform.
 
 \qua i)  If $u\in\rvo_j$, $u\ge 0$,  $B=\{P_h: u(P_h)>0\}$ and $n$
 is a positive integer $n$, then
 
 $$L^n_j(B)=\{P_h: T_{j;E;r}^n(u)(P_h)>0\}=\cf_n\text{-}L_{j^{(n)}}(B).$$

 \qua ii) for every $n\ge 1$, given
 $\{j,j'\}\in J$ we have $\{P_j,P_{j'}\}\in \whcg$
if and only if $L_j^n(P_{j'})\ne\empty$ if and only if there exists
$P_h\in\vo$ such that $T_{j;E;r}^n(\chi_{P_{j'}})(P_h)>0$. }

 \smad 
 i)  By (6.1) $T_{j;E;r}(u)(P_h)>0$ if and only if there exists $P_l\in B$ such that
 
 $T_{j;E;r} (\chi_{P_l})(P_h)>0$ if and only if $P_h\in L_j(B)$ by Lemma 6.1 i),
  and this proves the first equality if $n=1$.
 As $T_{j;E;r}^n(u)=\cf_n$-$T_{j^{(n)};E;r^n}(u)$ by Lemma 4.1 ii), the second equality
 follows from the first one with $n=1$ in $\cf_n$.
 The first equality for general $n$ follows by recursion. 
 Suppose the Lemma valid for some $n$. Then
 $T_{j;E;r}^{n+1}(u)(P_h)>0$ iff $T_{j;E;r}\big(T_{j;E;r}^n(u)\big)(P_h)>0$
 iff $P_h\in L_j(B_n)$ where $B_n=\{P_h: T_{j;E;r}^n(u)(P_h)>0\}$ using the case $n=1$.
 Thus, by the case $n$,  $T_{j;E;r}^{n+1}(u)(P_h)>0$ iff
 $h\in L_j\big(L_j^n(B)\big)=L_j^{n+1}(B)$.

\qua  ii) In view of Lemma 6.1 ii) in $\cf_n$ with $j^{(n)}$ in place of
$j$, this follows from i), putting $u=\chi_{P_j'}$ so that $B=\{P_{j'}\}$.
 \enpr

 \bigskip
  
  Let $\rvo_A:=\big\{u\in\rvo: u(P)=0\ \ \forall\, P\notin A\big\}$ when $A\subseteq \vo$.
  Clearly, $\rvo_A$ can be identified to $\bre^A$, and we will use such an identification.
  In particular, we will say that $u\in\rvo_A$ is {\it $A$-positive} if $u(P_h)>0$ for every
  $P_h\in A$.
 For every $s=1,...,m_j$, by Lemma 6.1 v) there exists a positive integer $n$
such that $L_j^{n}(C_{j,s})=C_{j,s}$.  Let $n_{j,s}$ 
be the minimum positive integer having such a property.

\bigskip
{\bf Lemma 6.3.}
{\sl  Suppose $E$ is an $r$-eigenform.

\qua i) We have $C_{j,s}=C'_{j,s}\cup C''_{j,s}$ where 
 
 $$C'_{j,s}=
 \Big\{P_{j'}\in C_{j,s}\:T_{j;E;r}^{n_{j,s}}(\chi_{P_{j'}})(P_h)>0\ \ \forall\,
 P_h\in C_{j,s}, T_{j;E;r}^{n_{j,s}}(\chi_{P_{j'}})(P_h)=0\ \  \forall\, 
 P_h\notin C_{j,s} \Big\},$$
  $$C''_{j,s}=\big\{P_{j'}\in C_{j,s}\:T_{j;E;r}^{n_{j,s}}(\chi_{P_{j'}})(P_h)=0\ \ \forall\,
 P_h\in \vo\big\}.$$
 
 Moreover, $C'_{j,s}\ne\empty$.
 
 \qua ii) If $\{j,j'\}\in J$,  then $\{P_j,P_{j'}\}\in\whcg$ if and only if $P_{j'}\in C'_{j,s}$.
 
 \qua iii) If $P_{j_1},P_{j_2}\in C'_{j,s}$ and $j_1\ne j_2$, then 
 $\{P_{j_1},P_{j_2}\}\in\whcg$.
 
 \qua iv) $T_{j;E;r}^{n_{j,s}}$ maps $\rvo_{C_{j,s}}$ into itself.
 
 \qua v) If $u\in \rvo$, $u\ge 0$ and $u(P)>0$ for at least one
 $P\in C'_{j,s}$, then $T_{j;E;r}^{n_{j,s}}(u)(P)>0$ for every $P\in C_{j,s}$. }
 \smad
 i) We have 
 
 $$C_{j,s}=L_j^{n_{j,s}}(C_{j,s})=\bigcup\limits_{P_{j'}\in C_{j,s}} 
 L_j^{n_{j,s}} (P_{j'}) \eqno (6.2)$$
 
 and, on the other hand, by Lemma 6.1 iv),  when
 $P_{j'}\in C_{j,s}$,  then $L_j^{n_{j,s}} (P_{j'})$  amounts either to $\empty$ or to
 $C_{j,s}$.
Using Lemma 6.2 i), with $u:=\chi_{P_{j'}}$, we see that in the former case
$P_{j'}\in C''_{j,s}$, and in the latter, $P_{j'}\in C'_{j,s}$. By (6.2),
the latter case occurs for at least one $P_{j'}\in C_{j,s}$.  
ii) This follows from i) and Lemma 6.2 ii).
iii) Given $j_1$, $j_2$ as in the statement iii), by Lemma 6.2 ii)

$$\cf_{n_{j,s}}\text{-}L_{j^{(n_{j,s})}}(P_{j_1})\cap
\cf_{n_{j,s}}\text{-}L_{j^{(n_{j,s})}}(P_{j_2})=
L_j^{n_{j,s}} (P_{j_1})\cap L_j^{n_{j,s}} (P_{j_2})=C_{j,s},$$

 thus,  iii) follows from  Lemma 6.1 iii) in $\cf_{n_{j,s}}$. iv) immediately follows from i),
 and v) follows from (6.1) and the definition of $C'_{j,s}$.
\enpr

 \bigskip
If $u\in\rvo$, $s=1,...,m_j$, let $g_{j,s}(u)=\big(u-u(P_j)\big) \chi_{C_{j,s}}
\in\rvo_{C_{j,s}}$
 Let $\tilde g_{j,s}: \rvo\to \rvo_{C'_{j,s}}$ be defined by
 $\tilde g_{j,s}(u)=u \chi_{C'_{j,s}}$.

 \bigskip
 {\bf Lemma 6.4} {\sl Suppose $E$ is an $r$-eigenform.
 
\qua  i) $\tilde g_{j,s}\circ T_{j;E;r}^{n_{j,s} }$ is a positive 
 linear operator from $\rvo_{C'_{j,s}}$ into itself, thus it has a unique 
 $C'_{j,s}$-positive eigenvector $\bar u_{j,s}$ of norm $1$ with eigenvalue $l_{j,s}>0$. 
 
 \qua ii) Putting $\tilde u_{j,s}:= T_{j;E;r}^{n_{j,s} }(\bar u_{j,s})$, $\tilde u_{j,s}$
 is a $C_{j,s}$-positive eigenvector of $T_{j;E;r}^{n_{j,s} }$ with eigenvalue
 $l_{j,s}$.
 
 \qua iii) for every $u\in \rvo_{C_{j,s}}$ 
 there exists $\pi_{j,s}(u)\in\bre$ such that
 ${T_{j;E;r}^{h n_{j,s} }(u)\over l_{j,s}^{h-1}}\llim\limits_
 {h\to +\infty} \pi_{j,s}(u)\tilde u_{j,s}.$}
 \smad
 i) Let $P_{j'}\in  C'_{j,s}$. By definition, we have 
 
 $$\tilde g_{j,s}\circ T_{j;E;r}^{n_{j,s} }(\chi_{P_{j'}})(P_h)=
 T_{j;E;r}^{n_{j,s} }(\chi_{P_{j'}})(P_h)>0\quad \forall\, P_h\in C'_{j,s}$$
 
 and i) is proved. ii) First,
 By Lemma 6.3 iv) and v), $\tilde u_{j,s}$ belongs to $\rvo_{C_{j,s}}$ and
 is $C_{j,s}$-positive. Next, note that 
 
 $$T_{j;E;r}^{n_{j,s} }\circ \tilde g_{j,s}=
 T_{j;E;r}^{n_{j,s} }\ \text{on}\ \rvo_{C_{j,s}}.\eqno (6.3)$$
 
  In fact,
 if $u\in \rvo_{C_{j,s}}$, then $u-\tilde g_{j,s}(u)=0$ 
 on $\vo\setminus\ C''_{j,s}$.
 Thus, we have for some $a_{j'}\in\bre$, 
  $u-\tilde g_{j,s}(u)=\sum\limits_{P_{j'}\in C''_{j,s}} a_{j'}\chi_{P_{j'}}$
  and thus, by definition of $C''_{j,s}$, 
  
  $$T_{j;E;r}^{n_{j,s} }\big(u-\tilde g_{j,s}(u)\big)=
  \sum\limits_{{P_{j'}}\in C''_{j,s}} a_{j'}T_{j;E;r}^{n_{j,s} }(\chi_{P_{j'}})=0.
  $$
 
and (6.3) is proved.  We have

$$T_{j;E;r}^{n_{j,s} } (\tilde u_{j,s})=T_{j;E;r}^{n_{j,s} }\big(
T_{j;E;r}^{n_{j,s} } (\bar u_{j,s})\big)=$$
$$T_{j;E;r}^{n_{j,s} }\big(\tilde g_{j,s}\circ
T_{j;E;r}^{n_{j,s} } (\bar u_{j,s})\big) 
=T_{j;E;r}^{n_{j,s} }(l_{j,s} \bar u_{j,s})=l_{j,s}\tilde u_{j,s}$$

and ii) is proved. iii) 
By the Perron-Frobenius Theory,
for every $u\in\rvo_{C'_{j,s}}$  there exists $\pi_{j,s}(u)$ such that 
 
 $${\big(\tilde g_{j,s}\circ T_{j;E;r}^{n_{j,s} })^h(u)\over l_{j,s}^h}\llim\limits_
 {h\to +\infty} \pi_{j,s}(u) \bar u_{j,s} \eqno (6.4)$$
 
More generally, (6.4) holds if $u\in \rvo_{C_{j,s}}$. In fact,  
$\tilde g_{j,s}\circ T_{j;E;r}^{n_{j,s} }(u)\in \rvo_{C'_{j,s}}$, thus it suffices to 
put $\pi_{j,s}(u)={1\over l_{j,s}}\pi_{j,s}\big(\tilde g_{j,s}\circ
T_{j;E;r}^{n_{j,s} }(u)\big)$.
 Let now $u\in \rvo_{C_{j,s}}$. By Lemma 6.3 iv) and (6.3) we have
 $T_{j;E;r}^{n_{j,s} }\Big(\big(\tilde g_{j,s}\circ T_{j;E;r}^{n_{j,s} })^h(u)\Big)=
T_{j;E;r}^{(h+1) n_{j,s} }(u)$.  Thus,

$${T_{j;E;r}^{(h+1) n_{j,s} }(u)\over l_{j,s}^h}
=T_{j;E;r}^{n_{j,s} }\bigg( {
\big(\tilde g_{j,s}\circ T_{j;E;r}^{n_{j,s} })^h(u)\over l_{j,s}^h}\bigg)\llim\limits_
 {h\to +\infty}\pi_{j,s}(u) T_{j;E;r}^{n_{j,s} }( \bar u_{j,s})= \pi_{j,s}(u) \tilde u_{j,s},$$
 
 and iii) is proved. \enpr

\bigskip
 {\bf Lemma 6.5.} {\sl Suppose  $E$ is an $r$-eigenform.
  
\qua i) We have $l_{j,s}=(\rho r_j^{-1})^{n_{j,s} }\in]0,1[$.
 
\qua  ii) If $X$ is a $T_{j;E;r}^{n_{j,s} }$-invariant  linear
  subspace of $\rvo_{C_{j,s}}$ and $u\in X$ and $P_j$ is an 
  $E$-nonharmonic point for $u$, then $\tilde u_{j,s}\in X$.}
\smad
By Lemma 6.4 ii), $T_{j;E;r}^{n_{j,s} }(\tilde u_{j,s})=l_{j,s} \tilde u_{j,s}$, thus,
by Lemma 2.6 i) we have
  
$$ L_{E}(\tilde u_{j,s})(P_j)=\big({r_j\over \rho}\big)^{n_{j,s} } 
L_E\Big(  T_{j;E;r}^{n_{j,s} }(\tilde u_{j,s})  \Big)(P_j)=$$
$$\big({r_j\over \rho}\big)^{n_{j,s} } 
   L_{E} (l_{j,s}\tilde u_{j,s})(P_j)=
\big({r_j\over \rho}\big)^{n_{j,s} } l_{j,s}  L_{E}(\tilde u_{j,s})(P_j). 
$$
   
On the other hand, by Lemma 6.4 ii)
$\tilde u_{j,s}$ is positive on the nonempty set
 $C_{j,s}$ and amounts to $0$ otherwise, thus, by Lemma 6.3 ii),
  we have $L_{E}(\tilde u_{j,s})(P_j)>0$, and i) follows at once.

\qua   I now prove ii). If $u\in X$, 
then by definition, also $\disp{{T_{j;E;r}^{h n_{j,s} }(u)\over l_{j,s}^{h-1}}}\in X$
for every positive integer $h$, and as $X$,
being a finite dimensional linear space, is closed, then by Lemma 6.4 iii)
$\pi_{j,s}(u)\tilde u_{j,s}\in X$, and in order to conclude it suffices 
to prove $\pi_{j,s}(u)\ne 0$.

Now, by  Lemma 6.4 iii) again, i) and Lemma 3.2, we have

$$\pi_{j,s}(u) L_{E}(\tilde u_{j,s}) (P_j)=L_{E}\big(\pi_{j,s}(u) \tilde u_{j,s}\big)(P_j)
=\lim\limits_{h\to+\infty}
L_{E}\Big( l_{j,s}^{-(h-1)}T_{j;E;r}^{h n_{j,s} }(u)\Big)(P_j)$$
$$=\lim\limits_{h\to+\infty}  l_{j,s}^{-(h-1)} 
L_{E}\Big( T_{j;E;r}^{h n_{j,s} }(u)\Big)(P_j)
=\big({\rho\over r_j}\big)^{n_{j,s} }L_E(u)(P_j)\ne 0$$

by the hypothesis that $P_j$ is an 
  $E$-nonharmonic point for $u$, thus $\pi_{j,s}(u)\ne 0$ and ii) is proved. \enpr

 \bigskip
 \centerline{\bf 7. Stable subsets of $\witi\cu$.}

 \medskip
 {\bf Lemma 7.1} {\sl If $E_1,E_2\in\td$ and $\cg_0(E_1)=\cg_0(E_2)=\whcg$, 
 and $u\in A^{\pm}(E_1,E_2)$,
  then $g_{j,s}(u)\in A^{\pm}(E_1,E_2)$
 for every  $s=1,...,m_j$.}
 \smad
 We have $E_i(u)=\sum\limits_{s=1}^{m_j} E_i\big(g_{j,s}(u)\big)$ for $i=1,2$. In fact,
 by the definition of $j$-component, we have $c_{j_1,j_2}(E_i)=0$ when
 $P_{j_1}$ and $P_{j_2}$ lie in different $j$-components. It follows
 
 $$E_i(u)=\sum\limits_{s=1}^{m_j}\bigg(\sum\limits_{
 P_{j_1}, P_{j_2}\in C_{j,s}} c_{j_1,j_2}(E_i) \big(
 u(P_{j_1})-u(P_{j_2})\big)^2+\sum\limits_{
 P_{j_1}\in C_{j,s}} c_{j_1,j}(E_i) \big(
 u(P_{j_1})-u(P_j)\big)^2\bigg)$$
 $$=\sum\limits_{s=1}^{m_j} E_i\big(g_{j,s}(u)\big),$$

 and the lemma easily follows. \enpr

 \bigskip
\qua Let $\witi \cu$ be the set of $(j,s)$ such that
$j\in\cu$ and $s=1,...,m_j$.
In general, for an arbitrary $r$-eigenform $E$,  
I say that a subset $B$ of $\witi\cu$ is $(E,r)$-stable if, for every $(j,s)\in B$
every $(j',s')\in\witi\cu\setminus B$,
 every $n\in\bna$ and every $i_1,...,i_n\in\cv$, 
 $P_{j'}$ is $E$-harmonic  for 
 $g_{j',s'} \big(T_{i_1,...,i_n;E;r}(\tilde u_{j,s})\big)$.

\bigskip
{\bf Lemma 7.2} {\sl Suppose $u\in \rvo$, and for every $(j,s)\in\witi\cu$,
 $P_j$ is $E$-harmonic point for 
 $g_{j,s}(u)$, with $E\in\td$. Then $u$ is constant on $\vo$. }
 \smad
 For every $(j,s)\in\witi\cu$, we have 
 
 $$0= L_E \big(g_{j,s}(u)\big)(P_j)=
  L_E\Big( \chi_{C_{j,s}}\big(u-u(P_j)\big) \Big)(P_j) $$
 
 so that, for every $j\in\cu$

$$L_E(u)(P_j)=L_E\big(u-u(P_j)\big)(P_j)=
L_E\Big(\sum\limits_{s=1}^{m_j} \chi_{C_{j,s}}\big(u-u(P_j)\big)\Big)(P_j)$$
$$=\sum\limits_{s=1}^{m_j} L_E\Big( \chi_{C_{j,s}}\big(u-u(P_j)\big) \Big)(P_j) =0$$

hence $u$ is $E$-harmonic at every point $P_j$, thus,
by Lemma 2.1, $u$ is constant on $\vo$. \enpr

\bigskip
{\bf Lemma 7.3} {\sl Suppose $E$ is an $r$-eigenform.
Then,  for every $(j,s)\in\witi\cu$, the function
$E_{j,s}$ defined by 

$$E_{j,s}(u)=\bigg(L_E\Big( T_{j;E;r}^{n_{j,s}}\big(g_{j,s}(u)\big)\Big)(P_j)\bigg)^2$$

has the form

$$E_{j,s}(u)=\sum\limits_{\{P_{j_i},P_{j_2}\}\in\whcg} 
d_{\{j_1,j_2\}} \big(u(P_{j_1})-u(P_{j_2})\big)^2$$

for suitable (possibly negative) $d_{\{j_1,j_2\}}\in\bre$.}
\smad
We have 
$g_{j,s}(u)= \disp{\sum\limits_{P_{j'}\in C_{j,s}} 
\big(u(P_{j'})-u(P_j)\big)\chi_{P_{j'}}}$, hence, 

$$T_{j;E;r}^{n_{j,s}} \big(g_{j,s}(u)\big)=\sum\limits_{P_{j'}\in C_{j,s}} 
\big(u(P_{j'})-u(P_j)\big)T_{j;E;r}^{n_{j,s}}(\chi_{P_{j'}} ),$$

and, as $T_{j;E;r}^{n_{j,s}}(\chi_{P_{j'}})(P_j) =\chi_{P_{j'}}(P_j)=0$, using Lemma 6.3 i),
we have

$$L_E\Big( T_{j;E;r}^{n_{j,s}}\big(g_{j,s}(u)\big)\Big)(P_j)=
\sum\limits_{h\in\cu\setminus\{j\}} \sum\limits_{P_{j'}\in C_{j,s}} 
 c_{j,h}(E) \big(u(P_{j'})-u(P_j)  \big)T_{j;E;r}^{n_{j,s}}(\chi_{P_{j'}})(P_h)$$
 $$=\sum\limits_{P_h\in C_{j,s}} \sum\limits_{P_{j'}\in C'_{j,s}} 
 c_{j,h} (E)\big(u(P_{j'})-u(P_j)  \big)T_{j;E;r}^{n_{j,s}}(\chi_{P_{j'}})(P_h)$$
 
 Thus, 
 
 $$E_{j,s}(u)=\sum\limits_{j'\in  C'_{j,s}} \alpha_{j'}
 \big(u(P_{j'})-u(P_j) \big)^2   +\sum\limits_{j_1,j_2\in   C'_{j,s}}
 \beta_{j_1,j_2} \big(u(P_{j_1})-u(P_j) \big)  \big(u(P_{j_2})-u(P_j) \big)$$
 
 for suitable $\alpha_{j'},  \beta_{j_1,j_2} \in\bre$. Now, using the identity
 $\disp{ab={1\over 2}\big(a^2+b^2-(a-b)^2\big)}$, in view of Lemma 6.3,
 ii) and iii), we conclude. \enpr

\bigskip
{\bf Theorem 7.4} {\sl If $E$ is an $r$-eigenform, then we have $r$-uniqueness
if and only if   every two   
 nonempty $(E,r)$-stable subsets of $\witi\cu$  
 are not disjoint.}
\smad
Suppose first every two   
 nonempty $(E,r)$-stable subsets of $\witi\cu$  
 are not disjoint, and prove the uniqueness.
 Let $E'$ be an $r$-eigenform.
Putting 

$$B^{\pm}:=\big\{(j,s)\in\witi\cu: \tilde u_{j,s}\in A^{\pm}(E,E')\big\},$$

 we prove that
$B^{\pm}$ are nonempty $(E,r)$-stable subsets of $\witi\cu$. They are nonempty. 
 In fact, let $u\in\witi A^{\pm}(E,E')$, and let $P_j$ be a minimum point of 
$u$. We can and do assume $u(P_j)=0$ and $u(P_{j'})>0$ 
with $\{P_j,P_{j'}\}\in\cg_0(E)$.
 In fact, given $\tilde j$ such that $u(P_{\tilde j})>0$
we can replace $u$ by $u-u(P_j)$, and consider a path connecting
$P_j$ with $P_{\tilde j}$ in $\cg_0(E)$ and consider two consecutive elements
$j'',j'$ in the path such that $u(P_{j''})=0$, $u(P_{j'})>0$, and finally replace $j''$ by $j$.
 Let $s=1,...,m_j$ be so that
$P_{j'}\in C_{j,s}$.  Then $P_j$ is not $E$-harmonic for
$g_{j,s}(u)$, and, in view of Lemma 7.1,
 $g_{j,s}(u)$ belongs to $A^{\pm}(E,E')\cap \rvo_{C_{j,s}}$, which is
a $T_{j;E;r}^{n_{j,s} }$-invariant  linear
  subspace of $\rvo_{C_{j,s}}$. Hence, by Lemma 6.5 ii),
   $\tilde u_{j,s}\in \witi A^{\pm}(E,E')$, and $(j,s)\in B^{\pm}$.

\qua We now prove that $B^{\pm}$ are  $(E,r)$-stable subsets of $\witi\cu$.
Suppose $(j,s)\in B^{\pm}$,
  $n\in\bna$,   $i_1,...,i_n\in\cv$, 
 and $P_{j'}$ is $(E,r)$-nonharmonic point for 
 
 $$\bar u:=g_{j',s'}\big(T_{i_1,...,i_n;E;r}(\tilde u_{j,s})\big),$$
 
and prove that  $(j',s')\in B^{\pm}$. By Prop.  2.2 i) and  Lemma 7.1, 
$\bar u\in A^{\pm}$, and also $\bar u\in \rvo_{C_{j',s'}}$. Moreover
$A^{\pm}(E,E')\cap \rvo_{C_{j',s'}}$,  is
a $T_{j';E;r}^{n_{j',s'}}$-invariant  linear
  subspaces of $\rvo_{C_{j',s'}}$. By Lemma 6.5 ii) again
  $\tilde u_{j',s'}\in \witi A^{\pm}(E,E')$,
  and $(j',s')\in B^{\pm}$.

\qua In conclusion, by the hypothesis, there exists
$(j,s)\in B^+\cap B^-$, thus $\tilde u_{j,s}\in \witi A^+(E,E')\cap \witi A^-(E,E')$, and
$E'$ is a multipole of $E$, thus we have $r$-uniqueness.

\qua Suppose conversely that there exist two
 nonempty $(E,r)$-stable subsets  $B_1$ and $B_2$ of $\witi\cu$ that are not disjoint,
 and prove that we have $r$-nonuniqueness. Let
 
 $$X_i:=\Big\{ u\in\rvo: 
 L_E\Big(g_{j,s} \big(T_{i_1,...,i_n;E;r}(u)\big)\Big)(P_j)=0\
 \forall\, i_1,...,i_n\in\cv,\ \forall\,  (j,s)\in\witi\cu\setminus B_i\Big\}.
 $$  

Then $X_1$ and $X_2$ are two $T_{E;r}$-invariant c-linear 
subspaces of $\rvo$. The only nontrivial point to prove is that they 
strictly contain the constants, but by definition of a $(E,r)$-stable subset, 
we have that $\tilde u_{j,s}\in X_i$ when $(j,s)\in B_i$. Moreover 
$X_1$ and $X_2$ are almost disjoint. In fact, if $u\in X_1\cap X_2$, 
then, as $B_1$ and $B_2$ are supposed to be disjoint, then 
in particular, $P_j$ is $E$-harmonic for
$g_{j,s}(u)$ for every $(j,s)\in\witi\cu$, thus, by Lemma 7.2,
$u$ is constant on $\vo$.

\qua By Lemma 2.4, the set $K_{t,E, X_2,X_1}$
 is closed and $\Lb_r$ invariant for every $t>1$. In view of Lemma 2.3,
 to prove
  the $r$-nonuniqueness, it suffices to prove that 
 for some $t>1$ it is also nonempty. Let
 
 $$E'=t\big(E-\delta E''\big), \quad E'':=\sum\limits_{(j,s)\in B_2}E_{j,s}.$$
 
We will prove that for suitable $\delta>0$ and $t>1$
 $E'\in K_{t,E, X_1,X_2}$. Note that for $\delta$ small enough, in view of Lemma
 7.3, $E'\in\td$.
 Next, by the definition of $E_{j,s}$ and Lemma 2.6 i), we have
 
 $$E''(u)=0 \iff   
P_ j\ \text{ is}\  E\text {-harmonic\ point\  for\  } g_{j,s}(u)
 \quad \forall\, (j,s)\in B_2. \eqno (7.1)$$
 
 In particular,  as $B_1$ and $B_2$ are disjoint,
 $E''=0$ on $X_1$, hence $E'=t E$ on $X_1$ is in any case satisfied, 
 and the condition  $E' \le t E$ on $\rvo$ is obviously satisfied.
 It remains to prove
 
 $$E'\le E\quad \text{on}\ X_2 \eqno (7.2)$$

Suppose $u\in X_2$, $u$ nonconstant. Then, by Lemma 7.2 and (7.1),
$E''(u)>0$, hence $(E-\delta E'')(u)<E(u)$. As the ratio
$\alpha:={E\over E-\delta E'' }$ 
has a minimum $t>1$ on the compact set
$S_{X_2}=:\{u\in X_2: u(P_1)=0, ||u||=1\}$, and
is $0$-homogeneous and invariant with respect 
to an additive constant, then  it has a minimum $t>1$
on the set of nonconstant functions of $X_2$. Clearly, with such a value of
$t$, (7.2) is satisfied.   \enpr

\bigskip
{\bf Remark 7.5}  By using the same argument as in Section 3, we need only finitely many 
calculations to verify whether a subset of $\witi\cu$ is $(E,r)$-stable.
In fact, we can see  that
 a subset $B$ of $\witi\cu$ is $(E,r)$-stable if and only if,
for every $(j,s)\in B$,
every $(j',s')\in\witi\cu\setminus B$,
 every $n=0,1,...,N-1$ and every $i_1,...,i_n\in\cv$, 
then $P_{j'}$ is $E$-harmonic  for 
 $g_{j',s'} \big(T_{i_1,...,i_n;E;r}(\tilde u_{j,s})\big)$.  \enpr

\bigskip
   
           \centerline {\bf References }
           \bigskip
         
           \rfbbkt\  Barlow M.T.,  Bass R.F., Kumagai T., and  Teplyaev A., 
           Uniqueness of Brownian Motion on Sierpinski Carpets. 
           J. Eur. Math. Soc. 12 (2010), 655-701.
               
           \smallskip
           \rfhmt\ Hambly,  B.M., Metz,  V.,  Teplyaev, A.: Self-Similar Energies on p.c.f. 
           Self-Similar Fractals. J. London Math. Soc. 74, 93-112 (2006).
           
           \smallskip
           \rfhh\  Hattori, K.,  Hattori, T., Watanabe, H.:  Gaussian Field Theories on 
           General Networks and the Spectral Dimension.  
           Progr. Theor. Phys. Suppl.  92, 108-143 (1987).
           
           \smallskip
           \rfki\  Kigami, J.:  Harmonic Calculus on p.c.f. Self-similar  Sets. 
           Trans. Amer. Math. Soc.  
           335,   721-755  (1993).
         
           \smallskip
           \rfkib\  Kigami, J.:  Analysis on fractals, Cambridge Tracts in Mathematics, 
           143. Cambridge University Press, Cambridge, 2001.
                    
          \smallskip
          \rfl\  Lindstr\o m, T.,  Brownian Motion on Nested Fractals. 
           Mem. Amer. Math. Soc.  No. 420
           (1990).
           
           \smallskip
           \rfm\  Metz, V.: How many Diffusions Exist on the Vicsek Snowflake?. 
            Acta Appl. Math.  32, 227-241  (1993).
                     
           \smallskip
           \rfmb\ Metz, V.: Renormalization Contracts on Nested Fractals.
           J. Reine Angew. Math.  480, 
              161-175  (1996).
              
          \smallskip
          \rfmc\  Metz, V.: The cone of diffusions on finitely ramified fractals. 
          Nonlinear Anal. 55  
              no. 6, 723--738 (2003).
              
            \smallskip
            \rfmd\ Metz, V.:  The short-cut test. Journal of Functional Analysis  220, 
            118-156 (2005).
            
           \smallskip
           \rfp\  Peirone, R.:  Convergence and Uniqueness Problems for Dirichlet Forms 
           on Fractals. Boll. U.M.I., 8, 3-B, 431-460  (2000).
          
          \smallskip
          \rfpb\ R. Peirone, Existence of Eigenforms on Nicely Separated Fractals, 
          in Proceedings of Symposia in Pure Mathematics,
          Amer. Math. Soc., Vol. 77, 231-241, 2008. 
      
           \smallskip
           \rfpc\  Peirone.: R.: Existence of Self-similar Energies on Finitely Ramified Fractals,
          to appear on Journal d'Analyse Mathematique.
          
          \smallskip
           \rfpd\ Peirone, R.: Uniqueness of Eigenforms on Fractals, Math. Nachr. 1-19 (2013),
           DOI 10.1002/mana.201200247
           
           \smallskip
           \rfsb\  Sabot, C.:  Existence and Uniqueness of Diffusions 
           on Finitely Ramified  
           Self-Similar Fractals. Ann. Sci. \'Ecole Norm. Sup. (4)  30 no. 5, 605-673  (1997).
           

           \smallskip
           \rfst\ Strichartz.: R. S.  Differential equations on fractals: A tutorial, 
           Princeton University Press, 2006.

 \end??????????????????????????????????????????????????

 i) is an immediate consequence of Lemma 2.3 i). By Lemma 3.2 we have
  
$$ L_{E}(\bar v_j)(P_j)={r_j\over \rho}    L_{E}(l_j\bar v_j)(P_j)=
{r_j\over \rho} l_j   L_{E}\bar v_j(P_j) 
$$
   
and thus, by i), $l_j=\rho r_j^{-1}$   and, in view also of Lemma 2.1,
 ii) is proved. I now prove iii). If $u\in X$, 
then by definition, also $\disp{T_{j,E;r}^k(u)\over l_j^k}\in X$, and as $X$,
being a finite dimensional linear space, is closed, then by Lemma 3.1
$\pi_j(u) \bar v_j\in X$, and in order to conclude it suffices 
to prove $\pi_j(u)\ne 0$.
Now, by  Lemma 3.1 again, ii) and Lemma 3.2, we have
$$\pi_j(u) L_{E}(\bar v_j) (P_j)=L_{E}\big(\pi_j(u) \bar v_j\big)(P_j)
=\lim\limits_{n\to+\infty}
L_{E}\Big( l_j^{-n}T_{j,E;r}^n(u)\Big)(P_j)$$
$$=\lim\limits_{n\to+\infty} l_j^{-n} L_{E}\Big(T_{j,E;r}^n(u)\Big)(P_j)
=L_E(u)(P_j)\ne 0$$

by the hypothesis that $P_j$ is an 
  $E$-nonharmonic point for $u$, thus $\pi_j(u)\ne 0$ and iii) is proved. \enpr

  \bigskip
  The idea of the use of Lemma 3.3 iii) is that if $X$ is $T_{E;r}$-invariant,
  then it contains $\bar v_j$ for every $P_j$ $E$-nonharmonic point
  of $T_{i,E;r}(u)$, then we can use the same process, 
  replacing $u$ by $\bar v_j$ and so on.
  In such a way, under suitable conditions a  nonempty $T_{E;r}$-invariant set contains
  so many $\bar v_j$ that two nonempty $T_{E;r}$-invariant sets
  have nonempty intersection, hence, by Theorem 2.8 we have uniqueness. 
  The details will be given in Section 4.
 
 \qua The well-known maximum principle states that $H_{1,E;r}(u)$ attains its
  maximum on $\von$ at points of $\vo$, and, further, when $E$ is positive,
  cannot attain its maximum at points of $\von\setminus\vo$ unless
  it is constant. Lemmas 3.4 and 3.5 can be seen as  
  variants of the maximum principle.

   \bigskip
  {\bf Lemma 3.4.} {\sl Let $u\in\rvo$, $r\in W$ and $v:=H_{1,E;r}(u)$. 
  Suppose $E\in\td$  is positive  and $Q\in \von\setminus\vo$.
  Suppose $S\subseteq st(Q)$ is such that
  
  \qua i) $Q$ is a maximum (resp. minimum) point for $v$
  on $V_i$ for every $i\in S$,
  
  \qua ii) $\psi_i^{-1}(Q)$ is $E$-harmonic for $T_{i,E;r}(u)$ for every $i\in st(Q)
  \setminus S$.
  
   Then $v=v(Q)$  on $V_i$ for every $i\in S$. }
  \smad   
   By (2.4), ii) and definition of $L_E$, putting $P_{j(i)}=\psi_i^{-1}(Q)$, we have
   
  $$\sum\limits_{i\in S}\sum\limits_{h\ne j(i)} r_i
  c_{j(i),h} \big(v(\psi_i(P_h))-v(\psi_i(P_{j(i)})\big)=0$$
  
  and the addenda are all nonpositive (resp. nonnegative), so that
  they all amounts to $0$, and as $E$ is positive,
then   $v(\psi_i(P_h))=v(Q)$ for all $P_h\in\vo$ and $i\in S$. \enpr

 \bigskip
 We will consider the graph $\cg_1$ on $\von$, formed by
 all sets of the form $\big\{\psi_i(P_{j_1}),\psi_i(P_{j_2})\big\}$ with $\{j_1,j_2\}\in J$,
 $i\in \cv$. 
 More generally, if $\cv'\subseteq\cv$, we consider the
 graph $\cg_1(\cv')$ on $V(\cv')$, formed by
 all sets of the form $\big\{\psi_i(P_{j_1}),\psi_i(P_{j_2})\big\}$ with $\{j_1,j_2\}\in J$,
 $i\in \cv'$. 
 Note that,
 if $Q,Q'\in V(\cv')$, $Q\ne Q'$, then $\{Q,Q'\}\in \cg_1(\cv')$ if and only if 
 $Q,Q'\in V_i$ for some $i\in\cv'$, if and only if $Q\in V(st(Q')\cap\cv')$.

 \qua Note that, if $\cv'$ is a connected subset of $\cv$, then it is easy to prove that
 $V(\cv')$ is  connected in $\cg_1(\cv')$. Moreover, we can take
 a path in $\cg_1(\cv')$ connecting two given points $Q,Q'\in V(\cv')$, with all points
 but the first and the last in $\von\setminus\vo$, using an argument like 
that in proof of Lemma 2.2.

\bigskip
{\bf Lemma 3.5.} {\sl Suppose $u\in\rvo$, $C$ is  a connected
subset of $\cv$,
$E$ is a positive $r$-eigenform, let $v:=H_{1,E;r}(u)$, and put 
$$ \partial_1=\big\{\psi_i(P_j)\in V(C)
:i\notin C,  j\  E\text{-nonharmonic\  for\ } T_{i,E;r}(u)\big\}\, (\subseteq
V(C)\setminus\vo),$$
$$\partial_2=\big\{P_j\in V(C), j\ E\text{-nonharmonic\ for\ } u\big\}
\, (\subseteq V(C)\cap\vo),$$
$$\partial=\partial_1\cup\partial_2 \, .$$

\qua  i) If $\partial_1\ne\empty$, and $v$ attains an
extremum (maximum or minimum) 
on $V(C)$ at points in $V(C)\setminus\partial_2$, then
 $v$ attains such an extremum at points in $\partial_1$.
 
\qua ii) If  $\partial_1=\empty$, and $v$ attains an
extremum (maximum or minimum) 
on $V(C)$ at points in $V(C)\setminus\partial_2$,
 then $v$ is constant on $V(C)$. In particular, if
 $\partial=\empty$, then $v$ is constant on $V(C)$.}
\smad
Suppose  $v$ attains its maximum 
or minimum $M$ on $V(C)$ at
 $\bar Q\in V(C)\setminus\partial$ and prove that

 $$v(Q)=v(\bar Q)\quad \forall\, Q\in V(st(\bar Q)\cap C).\eqno (3.1)$$
 
  In fact,
 if $\bar Q\notin\vo$, then, as $\bar Q\notin\partial_1$,  ii) of Lemma 3.4 holds with
 $\bar Q$ in place of $Q$ and (3.1) follows from Lemma 3.4,
 and  if on the contrary,
  $\bar Q=P_j\in\vo$, then $st(\bar Q)=\{j\}$, and as $\bar Q\notin\partial_2$,
 by Lemma 3.2 $j$ is $E$-harmonic for $v\circ\psi_j$, and 
 Lemma 2.3 i), $v\circ\psi_j$ is constant on $\vo$, that is
  $v=v(\bar Q)$ on $V_j=V\big(st(\bar Q)\cap C\big)$
and (3.1) holds. We can reformulate (3.1) in the following way.

$$v(Q)=M\quad\text{if}\  \{\bar Q,Q\}\in \cg_1(C).\eqno (3.2) $$

\qua Now,  let $\bar Q\in V(C)\setminus\partial_2$ 
be such that $v$ attains an extremum on $V(C)$ at $\bar Q$. Then,
we can connect 
$\bar Q$ to any point  $Q'\in V(C)$ by a path $(Q_0,...,Q_n)$
in $\cg_1(C)$ on $V(C)$ in such a way
that $Q_i\in V(C)\setminus\vo$ for $i=1,...,n-1$,
thus   $Q_i\in V(C)\setminus\partial_2$, for $i=0,...,n-1$.

\qua Hence, if $\partial_1=\empty$, $Q_i\in V(C)\setminus\partial$, 
for $i=0,...,n-1$, and, using (3.2), 
by recursion we have $v(Q')=v(Q_n)=v(\bar Q)=M$ and $v$ is constant 
on $V(C)$.
On the other hand, if  $\partial_1\ne\empty$, we take $Q'\in\partial_1$
and, in the path above we take $Q_l$ to be the first element
in $\partial_1$, possibly $l=0$,  thus, we consider the 
path $(Q_0,...,Q_l)$, and as $Q_i\in V(C)\setminus \partial$ for
$i=0,...,l-1$, we have $v(Q_l)=M$. \enpr

We define the function
$\bar v:\vin {n+1}\to\bre$   as

$$\bar v\Big(\psi_{i_1,...,i_n}(Q)\Big)=
H_{1;r;E}\big(H_{n;r^n;E}(u)\circ\psi_{i_1,...,i_n}\big)(Q)$$
for $Q\in \vin 1,\, i_1,...,i_n=1,...,k$, and prove that it
satisfies 

$$\bar v\in \cl(n+1,u), \quad S_{n+1,r^{n+1}}(E)(v)\ge  \Lb_{r^n}\big(\Lb_{r}(E)\big)(u)
\quad\forall v\in \cl(n+1,u)\, , \eqno(5.1)$$

and the equality holds if and only if $v=\bar v$.
 First, note that the definition of
$\bar v$ is correct, i.e., it does not depend on the representation of 
 $P\in\vin {n+1}$. Suppose 

$$P=\psi_{i_1,...,i_n}(Q)=\psi_{i'_1,...,i'_n}(Q')$$

 with $i_1,...,i_n, i'_1,...,i'_n=1,...,k,\,Q,Q'\in \vin 1$. 
We have either
 $(i_1,...,i_n)=(i'_1,...,i'_n)$, thus $Q=Q'$, or   by the nesting   axiom 
 $Q,Q'\in\vo$, and in any case 

$$H_{1;r;E}\big(H_{n;r^n;E}(u)\circ\psi_{i_1,...,i_n}\big)(Q)=
H_{1;r;E}\big(H_{n;r^n;E}(u)\circ\psi_{i'_1,...,i'_n}\big)(Q').$$

  Next, we see that $\bar v\in \cl(n+1,u)$.  Finally, if $v\in \cl(n+1,u)$, we have

$$S_{n+1;r^{n+1}}(E)(v)=
\sum\limits_{i_1,...,i_n=1}^k
r_{i_1}\cdots r_{ i_n}\Big(\sum \limits_{i_{n+1}=1}^k\
r_{i_{n+1}}E( v\circ\psi_{i_1,...,i_n}
\circ\psi_{i_{n+1}})\Big)$$
$$=\sum\limits_{i_1,...,i_n=1}^k r_{i_1}\cdots r_{ i_n} S_{1;r}(E)(v\circ\psi_{i_1,...,i_n})
\ge \sum\limits_{i_1,...,i_n=1}^k r_{i_1}\cdots r_{ i_n}\Lb_{r}(E) 
\big(v_{|\vin n}\circ\psi_{i_1,...,i_n}\big)$$
$$=S_{n;r^n}(\Lb_{r}(E))\big(v_{|\vin n}\big)\ge 
\Lb_{r^n}\big(\Lb_{r}(E)\big)(u)\, .$$
and (5.1) holds. Moreover, the first inequality is in fact an equality if and only if 

$$v\circ \psi_{i_1,...,i_n}=H_{1;r;E}
\big({v_{|\rvin n}\circ\psi_{i_1,...,i_n}}\big)$$
on $\vin 1$ for all $i_1,...,i_n=1,...,k$, the second is an equality
 if and only if $v_{|\vin n}=H_{n;r^n;E}(u)$, if and only if

$$v\circ \psi_{i_1,...,i_n}=H_{n;r^n;E}(u)\circ \psi_{i_1,...,i_n}$$
on $\vo$ for all $i_1,...,i_n=1,...,k$. Hence, 
the equality holds in (5.1) if 
and only if $v=\bar v$.  As a consequence, by the definition of $\Lb_{r^{n+1}}$ we have

$$\Lb_{r^{n+1}}(E)(u)=S_{n+1;r^{n+1}}(E)(\bar v)=\Lb_{r^n}\big(\Lb_{r}(E)\big)(u)$$

and i) follows. On the other hand, $\bar v=H_{n+1;r^{n+1};E}(u)$, hence, using 
the definition of $\bar v$ at the point $Q=\psi_{i_{n+1}}(P_j)$, we get

$$H_{n+1;r^{n+1};E}(u)\circ\psi_{i_1,...,i_n,i_{n+1}}(P_j)=$$
$$\bar v\circ \psi_{i_1,...,i_n,i_{n+1}}(P_j)=T_{i_{n+1};E;r}\big(
H_{n;r^n;E}(u)\circ\psi_{i_1,...,i_n}\big) (P_j)$$

and ii) follows by recursion.\enpr

I and
 let $(\bar{re})_{j,s}(u)=u|_{C_{j,s}}$.  
 If $u\in\bre^{C_{j,s}}$ let $\iota_{j,s}(u)\in\rvo$ be defined by
 $\iota_{j,s}(u)(P_h)=\cases u(P_h)&\ \text{if}\  P_h\in C_{j,s}\\
 0&\ \text{otherwise} \endcases.$ Note that
  
 $$g_{j,s}=\iota_{j,s}\circ (\bar{re})_{j,s}.$$

\bigskip
Now,
 if $L_j^{n}(C_{E;j,s})\supseteq C_{E;j,s}$, then
 $\zeta_{j,s}(T_{j;E,r}^{n})$ is a positive linear operator from
$\bre^{C_{E;j,s}}$ into itself. I will say that $C_{E,j,s}$ is {\it regular} if
or a regular $j$-component if there exists a positive integer $n$ such that
$L_j^{n}(C_{E;j,s})\supseteq C_{E;j,s}$. In such a case I will say that
$n$ is a regularity index of $(j,s;E,r)$.

 \bigskip
 {\bf Lemma 7.4} {\sl  Let $E\in\td$. 
 Then, either   $L_j^{m_j}(\cu\setminus\{j\})\cap C_{E;j,s}=\empty$ or there exists
  $s'$ (that can be equal to or different from $s$) such that
  
  $$L_j^{m_j}(C_{E;j,s'})\supseteq C_{E;j,s}, \quad
L_j^{n}(C_{E;j,s'})\supseteq C_{E;j,s'} {\ for\ some\ } n=1,...,m_j. \eqno(6.1)$$}
  \smad
 If $L_j^{m_j}(\cu\setminus\{j\})\cap C_{E;j,s}\ne \empty$, then there exists
 $s''$ such that $L_j^{m_j}(C_{E;j,s''})\supseteq C_{E;j,s}$,
 Then there exist $s'_0=s'', s'_1,s'_2,...,s'_{m_j}=s$ such that
 $L_j (C_{E;j,s'_l})\supseteq C_{E;j,s'_{l+1}}.$
 Then, let $s'_h=s'_k$ with $0\le h<k\le m_j$. 
 We put $s_l=s'_l$ for $l=k,...,m_j$, and
  $s_{h-1}=s'_{k-1}$,
 and more generally, $s_{h-l(k-h)-r}=s'_{k-r}$ when
 $l=0,1,2....$, $r=0,...,k-h-1$, and $s'=s_0$. \enpr

 \bigskip
 Given $L:\rvo\to\rvo$, let $\zeta_{j,s}(L)=(\bar{re})_{j,s}\circ L\circ \iota_{j,s}$
 from $\bre^{C_{E;j,s}}$ into itself.

 \bigskip
 {\bf Lemma 7.5.} {\sl If $L_j^{n_i}(C_{E;j,s})\supseteq C_{E;j,s}$, $i=1,2$, then 

$$\zeta_{j,s}\big(T_{j;E;r}^{n_1})\circ \zeta_{j,s}(T_{j;E;r}^{n_2}\big)=
 \zeta_{j,s}\big(T_{j;E;r}^{n_1+n_2}\big). $$}
 
 %
 
 %
 \smad
 We prove that, if $u\in\rvo$ and $u(P_h)=0$ if $P_h\in C_{E;j,s}$, then
 
 $$T_{j;E;r}^{n_1}(u)(P_h)=0 \quad\forall\, P_h\in C_{E;j,s}.\eqno (7.2)$$
 
  In fact, $u=\sum\limits_{l\notin C_{E;j,s}} u(P_l) e_l$, thus
  
  $$T_{j;E;r}^{n_1}(u)(P_h)=\sum\limits_{l\notin C_{E;j,s}} u(P_l)
T_{j;E;r}^{n_1} (e_l)(P_h)$$

and $T_{j;E;r}^{n_1} (e_l)(P_h)=0$ as, if not, by Lemma 6.2, $h\in L_j^{n_1}(j_1)$,
for some $j_1\notin C_{E;j,s}$. Hence, by Lemma 6.1 iv), 
$L_j^{n_1}(j_1)\supseteq C_{E;j,s}$, but also, by hypothesis, 
$L_j^{n_1}(j_2)\supseteq C_{E;j,s}$ for some $j_2\in C_{E;j,s}$, so that,
by Lemma 6.1 v),  $\{P_{j_1},P_{j_2}\}\in
\Lb\big(\cg_0(E)\big)$, hence, as $j_2\in C_{E;j,s}$, also
$j_1\in C_{E;j,s}$ , a contradiction, and (6.2) is proved. Now,

$$\zeta_{j,s}(T_{j;E;r}^{n_1})\circ \zeta_{j,s}(T_{j;E;r}^{n_2})(u)-
 \zeta_{j,s}\big(T_{j;E;r}^{n_1+n_2}\big)(u)=(\bar{re})_{j,s}\circ
T_{j;E;r}^{n_1}\Big((g_{j,s}-Id)\circ  T_{j;E;r}^{n_2}\circ\iota_{j,s}(u)\Big)$$
 
 But, by definition $g_{j,s}-Id=0$ on $C_{E;j,s}$, hence, by (7.2)
 $T_{j;E;r}^{n_1}\Big((g_{j,s}-Id)\circ  T_{j;E;r}^{n_2}\circ\iota_{j,s}(u)\Big)=0$
 on $C_{j,s}$, thus  $(\bar{re})_{j,s}\circ
T_{j;E;r}^{n_1}\Big((g_{j,s}-Id)\circ  T_{j;E;r}^{n_2}\circ\iota_{j,s}(u)\Big)=0$. \enpr

 \bigskip
 {\bf Lemma 7.6} {\sl  For every $j\in\cu$ and every $s=1,...,m_j$ we have
 
 \qua i)  $L_j(C_{E;j,s})\ne\empty$ for every 
 $s=1,...,m_j$
 
 \qua ii) For every $j\in \cu$, $s=1,...,m_j$, 
 there exists a positive integer $n$ such that
 $L_j^n(C_{E;s,j})$ contains a regular  $j$-component.}
 \smad
 By definition of  component there exists $j'\in C_{E;s,j}$ such that
 $\{j,j'\}\in \cg_0 (E)=\Lb\big(\cg_0(E)\big)$. Now, i)
 follows from Lemma 7.1 ii).
 
 ii) Let $s_0=s$ and $s_h$ be so that  
 $L_j(C_{E;s_h,j})\supseteq C_{E;s_{h+1},j}$. We have $s_h=s_k$ for some $0\le h<k$.
 Then $L_j^h(C_{E;s,j})\supseteq C_{E;s_h,j}$ and 
 $L_j^{h-k}(C_{E;s_h,j})\supseteq C_{E;s_h,j}$, so that $C_{E;s_h,j}$ is a regular
 $j$-component. \enpr

\bigskip
{\bf Lemma 7.7} {\sl If $C_{E;j,s}$ is regular, then there exists a unique
$\bar u_{j,s}$ which is a positive eigenvector of norm $1$
of $\zeta_{j,s}(T_{j;E,r}^{n})$ with $n$ regularity index of 
$(j,s;E,r)$.}
\smad
Given $n$ regularity index of $(j,s;E,r)$, then $\zeta_{j,s}(T_{j;E,r}^{n})$
has a unique positive eigenvector of norm $1$. It remains to prove that 
such an eigenvector does not depend on $n$. Suppose $n_1, n_2$ regularity indices
of $(j,s;E,r)$ and $u^{(n_i)}$ is 
the positive eigenvector of norm $1$ of $\zeta_{j,s}(T_{j;E,r}^{n_i})$. We are going
 to prove that $u^{(n_1)}=u^{(n_2)}$. It suffices to note that,
 in view of Lemma 7.5,  both 
 $u^{(n_1)}$ and $u^{(n_2)}$ are positive eigenvectors
 of norm $1$ of $\zeta_{j,s}(T_{j;E,r}^{n_1 n_2})$, and clearly
 $n_1 n_2$ is a regularity index of $(j,s;E,r)$. \enpr

\bigskip
\qua Let $\tilde u_{j,s}=\iota(\bar u_{j,s})$, and let $\witi \cu$ be the set of $(j,s)$ such that
$j\in\cu$ and $C_{E;j,s}$ is regular.
In general, for an arbitrary $r$-eigenform $E$,  
I say that a subset $B$ of $\witi\cu$ is $(E,r)$-stable if, for every $(j,s)\in B$,
 every $n\in\bna$ and every $i_1,...,i_n\in\cv$, 
 when $j'$ is $E$-nonharmonic point for $g_{j',s'}\Big(T_{i_1,...,i_n,E;r}(\tilde u_{j,s})\Big)$
and $C_{E;j',s'}$ is regular, then $(j',s')$  belongs to $B$.

\smad
By definition of $(E,r)$-hstable,
we have to prove that if $B$ is hstable, then  if $u\in\rvo$
 is such that
 every point of $\cu\setminus B$ is $E$-harmonic for $u$, then
 every point of $\cu\setminus B$ is $E$-harmonic for $T_{i,E;r}(u)$ as well.

\qua  Suppose  $C$ is a component of $\cv\setminus A$ and
  use Lemma 3.5 and notation there.
   We have $\partial_2=\empty$ as, if $P_j\in V(C)$, then
 $j\in C$, thus  $j\notin\cu\setminus B$, and $j$ is $E$-harmonic for $u$.
 Also, if $Q\in\partial_1$, as $C$ is a component of
 $\cv\setminus A$ and by the definition of $\partial_1$, then
 $Q\in V(C)\cap V(A)$. Thus, by i) $\partial_1$ has at most one element.
 By Lemma 3.5 $H_{1,E,r}(u)$ is constant on $V(C)$, as if
 $\partial_1\ne\empty$, then the maximum and the minimum on
 $V(C)$ are attained at the same point.
   Therefore, if $i\in C$ then
  $T_{i,E;r}(u)$ is constant, thus
  
  $$j  \ \text{ is}\  E\text{-harmonic \ for\ } T_{i,E;r}(u)\ \ \forall j=1,...,N, 
i \in\cv\setminus   A\, . \eqno (4.1)$$
  
   It remains to prove that,
  if $\bar i\in A$ and $j\in \cu\setminus B$, then $j$ is
  $E$-harmonic for $T_{\bar i,E;r}(u)$. 
  If $\psi_{\bar i}(P_j)\in\vo$, then $\bar i=j\in\cu\setminus  B$,
  thus $j$ is $E$-harmonic for $u$, hence, by  Lemma 3.2,
  also for $T_{\bar i,E;r}(u)$. 
  Suppose on the contrary, $\psi_{\bar i}(P_j)\notin\vo$.
   Then, by ii) in the definition of hstable,
   if $i\in st\big(\psi_{\bar i}(P_j)\big)\setminus\{\bar i\}$ 
   then $i\notin A$, thus, by (4.1),  
  $L_E\big(T_{i,E;r}(u)\big)\big(\psi_i^{-1}(\psi_{\bar i}(P_j)\big)=0$.
Hence by (2.4)   we have $L_E\big(T_{\bar i,E;r}(u))(P_j)=0$,
 and we conclude. \enpr

  \bigskip
  {\bf Remark 4.3.} In  Lemma 4.2, in order to prove that $H_{1,E;r}(u)$
   is constant on $V(C)$, we never use ii) in definition of a hstable set. Thus,
   if $B$ and $A$ satisfy i) in the definition of a hstable set  and
   $u\in \bu_{B,E}$, then  $H_{1,E;r}(u)$ is constant on 
   $V(C)$ for every $C$ component of
   $\cv\setminus A$. \enpr

\bigskip
 {\bf Lemma 4.5.} {\sl If $E$ is a positive $r$-eigenform and there exist two  disjoint
 $(E,r)$-hstable subsets $B$ and $B'$ of $\vo$  having at least two elements,
 then there is $r$-nonuniqueness.}
 \smad
 The sets  $\bu_{B,E}$ and $\bu_{B',E}$
 are $T_{E;r}$-invariant by definition, and they are trivially linear subspaces. Moreover,
 they strictly contain the constants, as if $j_1, j_2$ are two different
 elements of $B$, then there exists
 a function $u$ taking arbitrary (that hence we can and do assume different) values
 at $P_{j_1}$ and $P_{j_2}$ such that every $P_j$ with $j_1\ne j\ne j_2$
 is $E$-harmonic for $u$.
 This is a well-known fact: we can take the function
 on $\vo$ taking the prescribed values at $P_{j_1}$ and $P_{j_2}$ minimizing
 $E$. Thus $u\in\bu_{B,E}$ and $u$ is nonconstant,
  and  a similar argument works for $\bu_{B',E}$. Finally,
  $\bu_{B,E}$ and $\bu_{B',E}$ are almost disjoint, as,
  if both every point of $\vo\setminus B$ and every point of
  $\vo\setminus B'$ is $E$-harmonic for $u$, then every point of $\vo$ 
  is $E$-harmonic for $u$, since $B$ and $B'$ are disjoint, hence $u$ is constant.
  We now immediately conclude using Theorem 2.8. \enpr

 \bigskip
 The following corollary is specially useful as it allows us to establish
 $r$-nonuniqueness for every $r$ for which we have $r$-existence,
 without knowing an $r$-eigenform.

 \bigskip
  {\bf Corollary 4.6.} {\sl If   (A) holds and there exist two  disjoint
 hstable subsets of $\vo$  having at least two elements, then 
   $r$-existence implies $r$-nonuniqueness.}
 \smad
 This immediately follows from  Lemma 4.2 and Lemma 4.5. \enpr

 \bigskip
 As we will see in Theorem 5.3, in some specific cases condition 
 of Corollary 4.6 is necessary and sufficient to have nonuniqueness. However,
 although I suspect that such situations could be relatively frequent, I do not see
 at this moment how to characterize them.

\qua I say that a pair $\{j_1,j_2\}$ with $j_1,j_2=1,...,N$, $j_1\ne j_2$ is a diagonal if 
there exists a path 
$(i_0=j_1,i_1,...,i_n=j_2)$ such that 

\qua i) $V_{i_h}\cap V_{i_{h+1}}\ne\empty$ if $h=0,..,n-1$

\qua ii) If $\psi_{i_h}(P_j)\in V_{i_{h'}}$ with $h'\ne h$, then $j\in\{j_1,j_2\}$.   
 
 In such a case, I say that $\{j_1,j_2\}$ is a proper diagonal if, putting
 $B=\{j_1,j_2\}$,  $A=\{ i_0, i_1,...,i_n\}$, then i)  of definition of hstable set
 holds.  I say that 
  $\{j_1,j_2\}$ is {\it weakly proper} if,
  whenever $C$ is a component of $\cv\setminus A$ then
   there exists $i\in A$ such that $V(C)\cap V(A)\subseteq V_i$.
 Of course, a proper diagonal is weakly proper. The converse holds
 provided
 
 $$C \text{\ component\ of\ } \cv\setminus A \Rightarrow
 \#\big(V(C)\cap V_i\big)\le 1\quad \forall\, i\in A\eqno (K)$$
 
 Condition (K) is satisfied relatively frequently.

 \bigskip
 {\bf Corollary 4.7.} {\sl If (A) holds and
 there are two disjoint proper diagonals, then for every $r$ 
  $r$-existence implies $r$-nonuniqueness.} 
 \smad 
 This is a particular case of Corollary 4.6, as ii) in definition of hstable set
 follows from ii) in definition of a diagonal. \enpr

 \bigskip 
\qua  Note that in view of Corollary 4.7 we can find fractals
 with $r$-nonuniqueness which are very far from being a tree,
 and this is specially true for $N$ large. In fact,
 we can have two proper diagonals with at the same time no constraints
 on the vertices not on the diagonals.
 \bigskip

           \centerline {\bf References }
           \bigskip
         
           \rfbbkt\  Barlow M.T.,  Bass R.F., Kumagai T., and  Teplyaev A., 
           Uniqueness of Brownian Motion on Sierpinski Carpets. 
           J. Eur. Math. Soc. 12 (2010), 655-701.
         
           \rfbft\ Boyle, B.,  Cekala,  K., Ferrone, D.,   Rifkin, N.,  Teplyaev, A.: 
           Electrical Resistance of N-gasket Fractal Networks.   
           Pacific Journal of Mathematics, 233, 
           No 1, 15-40 (2007).
          
           \smallskip
           \rfhmt\ Hambly,  B.M., Metz,  V.,  Teplyaev, A.: Self-Similar Energies on p.c.f. 
           Self-Similar Fractals. J. London Math. Soc. 74, 93-112 (2006).
           
           \smallskip
           \rfhh\  Hattori, K.,  Hattori, T., Watanabe, H.:  Gaussian Field Theories on 
           General Networks and the Spectral Dimension.  
           Progr. Theor. Phys. Suppl.  92, 108-143 (1987).
           
           \smallskip
           \rfki\  Kigami, J.:  Harmonic Calculus on p.c.f. Self-similar  Sets. 
           Trans. Amer. Math. Soc.  
           335,   721-755  (1993).
         
           \smallskip
           \rfkib\  Kigami, J.:  Analysis on fractals, Cambridge Tracts in Mathematics, 
           143. Cambridge University Press, Cambridge, 2001.
                    
          \smallskip
          \rfl\  Lindstr\o m, T.,  Brownian Motion on Nested Fractals. 
           Mem. Amer. Math. Soc.  No. 420
           (1990).
           
           \smallskip
           \rfm\  Metz, V.: How many Diffusions Exist on the Vicsek Snowflake?. 
            Acta Appl. Math.  32, 227-241  (1993).
                     
           \smallskip
           \rfmb\ Metz, V.: Renormalization Contracts on Nested Fractals.
           J. Reine Angew. Math.  480, 
              161-175  (1996).
              
          \smallskip
          \rfmc\  Metz, V.: The cone of diffusions on finitely ramified fractals. 
          Nonlinear Anal. 55  
              no. 6, 723--738 (2003).
              
            \smallskip
            \rfmd\ Metz, V.:  The short-cut test. Journal of Functional Analysis  220, 
            118-156 (2005).
            
           \smallskip
           \rfp\  Peirone, R.:  Convergence and Uniqueness Problems for Dirichlet Forms 
           on Fractals. Boll. U.M.I., 8, 3-B, 431-460  (2000).
          
          \smallskip
          \rfpb\ R. Peirone, {\it Existence of Eigenforms on Nicely Separated Fractals}, 
          in Proceedings of Symposia in Pure Mathematics,
          Amer. Math. Soc., Vol. 77, 231-241, 2008. 
      
           \smallskip
           \rfpc\  Peirone, R.: Existence of Self-similar Energies on Finitely Ramified Fractals,
           preprint at hal.archives-ouvertes.fr/docs/00/62/86/61/PDF/exist4.pdf, 
          October 2011.

           \smallskip
           \rfs\  Sabot, C.:  Existence and Uniqueness of Diffusions 
           on Finitely Ramified  
           Self-Similar Fractals. Ann. Sci. \'Ecole Norm. Sup. (4)  30 no. 5, 605-673  (1997).
           
            \smallskip
          \rfsb\  Sabot, C.: Espaces de Dirichlet reli\'es par des points 
         et application aux 
          diffusions sur les fractals finiment ramifi\'es. Potential Anal. 11 no. 2, 
          183--212  (1999).

           \smallskip
           \rfst\ Strichartz, R.: Isoperimetric estimates on Sierpinski gasket type fractals.
           Trans. Amer. 
           Math. Soc. 351  no. 5, 1705--1752 (1999). }

  \end
   
 {\bf Theorem 3.1.} {\sl Let $\cf$ be a generalized Vicsek set and let $r\in W$.
 If it is a quasi-tree and
  $\cv$ contains$D^+$ and $D^-$, and  there exists a $\Lb_r$-eigenform, 
  then  there exist infinitely
  many nonproportional $\Lb_r$-eigenforms, in particular there exist 
  infinitely many nonproportional
  eigenforms.  Conversely,if  we  have two nonproportional $\Lb_r$-eigenforms, 
  then it is a quasi-
  tree and $\cv$  contains $D^+$ and $D^-$.}
 
 \smallskip
  Proof. If $\cv=\widehat {\cv}:=D^+\cup D^-$, I consider the effective resistances 
  $R_{j_1,j_2}$, and use arguments in \rfsb.Let $\bar R$ be the map from $[0,+\infty[^6$ 
  to $\bre^6$ that sends the conductances $c_{j_1,j_2}$ to
   the effective resistances $R_{j_1,j_2}$. 
  If $\cv$ is a tree,  in view of  \rfsb, $c=(c_{j_1,j_2})$ is a $\Lb_r$-eigenform if and only if 
  $\bar R(c)$  satisfies a homogeneous linear system $(S)$. Moreover,  the equations
   corresponding to the diagonals are of the forms $\alpha R_{1,3}=R_{1,3}$, 
   $\beta R_{2,4}=R_{2,4}$. Now, if there exists a $\Lb_r$-eigenform with  coefficients 
   $(\bar c_{j_1,j_2})$, then  $(S)$ has a solution $\bar R(\bar c)$, 
   so that $\alpha=\beta=1$. 
   Therefore, $(S)$ has a two-dimensional subspace of solutions. 
   Since $\bar R$ is continuous 
   (trivial) and one-to-one (by \rfkib , Theorem 2.1.12), then
    it is open in $]0,+\infty[^6$ by a well 
   known result in euclidean topology. As moreover, $\bar R$ 
   is $(-1)$-homogeneous, it follows
    that we have infinitely many mutually independent 
    $\Lb_r$-eigenforms. If $r_i=1$ for each $i$,
  then by the symmetry it is easy to see that the form 
  with all coefficients $c_{j_1,j_2}=1$ is an
 eigenform. If, more generally, $\cf$ is a quasi-tree 
 and $\cv$ contains$D^+$ and $D^-$, then we
  easily see that in the definition of a quasi-tree, we can take
   $\widetilde {\cv}=\widehat {\cv}$, 
  and the operator $\Lb_r$ is the same with respect to $\widetilde {\cv}$ 
  and to $\widehat {\cv}$, 
  so that we have the same eigenforms as before.  
  
  \quad$\,$ To prove the converse,
  suppose $E,E'$ are two nonproportional eigenforms,  put $T_i=T_{i;E;r}$, and let 
  $u\in \ta^+(E,E')$. We can assume $u$ takes its minimum at $P_1$. 
  Lemma 2.2 provides 
  (at least) two elements of $\ta^+(E,E')$, namely $u_j$ for $P_j$ 
  extremum points of $u_1$. 
  Similarly, we find two elements of $\ta^-(E,E')$ of the form $u_j$. Since $T_1(u_1)$ is, 
  by definition, a positive multiple of $u_1$, and by Lemma 2.3, as $\psi_1(P_2)$ 
  and $\psi_1(P_4)$ are $1$-points, $T_1(u_1)$ cannot take its extrema  at $P_2$
   and $P_4$, 
  then $u_1,u_3\in \ta^+(E,E')$. On the other hand, using Lemma 2.1 iv and Lemma 2.2, 
  $\ta^+(E,E')$ contains every $u_j$ with $P_j$ extremum point of  
  nonconstant $T_i(u_1)$. 
  As we only have four $u_j$, in view of Prop. 2.1 iii, then every nonconstant $T_i(u_1)$ 
  takes its extrema  at $P_1$ and $P_3$.
  
  \quad$\,$  Next, $v_1:=H_{1;E;r}(u_1)$ is constant on any increasing diagonal
   $D\ne D_+$.  In fact, let $D=D_{+,h,a,b}$. Since $D\ne D^+$, 
  at least one of the two points  $Q:=\psi_{(a,a+h)}(P_1)$,
  $Q':=\psi_{(b,b+h)}(P_3)$ is not in $\vo$. 
  Suppose for example $Q\in\von\setminus\vo$. Since $v_1$  
  attains its maximum or minimum 
  on $V_{(a,a+h)}$ at $Q$ and, clearly, $Q$ is a $1$-point, 
  by Lemma 2.3 $v_1$ is constant 
  on  $V_{(a,a+h)}$. By the same argument we can now deduce that $v_1$ is constant on 
   $V_{(a+1,a+1+h)}$ (if $a+1\le b$), and iterating the argument, 
   $v_1$ is constant on $D$, 
    as claimed.\quad$\,$ It follows in particular that $\cv$  contains $D_+$, 
    since otherwise, 
    $v_1$ would be constant on every  increasing diagonal, thus, 
    by a connectedness argument,  
    on $\von$, which is impossible. Note that since $v_1$ 
    attains its minimum on $V_{(0,0)}$ 
    at $\psi_{(0,0)}(P_1)=P_1$, it attains its maximum 
    on $V_{(0,0)}$ at $\psi_{(0,0)}(P_3)$, 
    and  by Lemma 2.3 again,  on $V_{(a,a)}$, $v_1$ is not 
    constant and attains its minimum  
    at $\psi_{(a,a)}(P_1)$ and its maximum at $\psi_{(a,a)}(P_3)$, for each $a=0,...,2n$. 
    Since $\psi_{(a,a)}(P_3)=\psi_{(a+1,a+1)}(P_1)$ we thus have

    $$v_1(\psi_{(a+1,a+1)}(P_j))>v_1(\psi_{(a,a)}(P_j)) \eqno (3.1)$$
    
    for every given $j=1,2,3,4$. Now, if there exists a path connecting 
    two different $1$-cells in $D^+$, 
    disjoint from $D^+$, then its extension  contains two different points either of the form 
    $\psi_{(a,a)}(P_2)$,  or of the form $\psi_{(a,a)}(P_4)$, and $v_1$
     takes the same value 
    on such two points, contrary to (3.1). Thus, I have proved that 
    there exists no path connecting 
    two different $1$-cells in $D^+$, disjoint from $D^+$. 
    
    \quad$\,$ Similarly,
    we can prove  that $\cv$ 
    contains $D^-$ and that there exists no path connecting two different 
    $1$-cells in $D^-$, 
    disjoint from $D^-$. It easily follows that $\cv$ is a quasi-tree. \enpr

    \medskip
    {\bf Remark 3.2.} The argument in Theorem 3.1 shows, more generally, 
    that we have the nonuniqueness  of the $\Lb_r$-eigenformon every quasi-tree 
    fractal containing two different long diagonals and having a $\Lb_r$-eigenform 
    with positive coefficients. Here, by long ($j_1,j_2$) diagonal I mean a sequence 
   $({i_1},...,{i_h})$ with $i_1=j_1$, $i_h=j_2$, $\psi_{i_l}(P_{j_2})=\psi_{i_{l+1}}(P_{j_1})$.
    In fact, I never used in that part of the proof the fact that the fractal 
    is a generalized Vicsek set,  neither that $N=4$. \enpr 
     
     \bigskip

     \bigskip
Now, we can give some useful information about $E$-nonharmonic points.

  \bigskip
  {\bf Lemma 3.6.} {\sl Suppose  $u\in\rvo$.

 \qua  i) If $E\in\td$ and $Q\in V_i\setminus\vo$ is $1$-point, then 
 $\psi_i^{-1}(Q)$ is an $E$-harmonic  point for $T_{i,E;r}(u)$.

\qua  ii) if $E$ is a positive
  $r$-eigenform, and $Q\in V_i\cap V_{i'}$, $i'\ne i$, 
and $Q$ is an extremum point 
 for $H_{1,E;r}(u)$ on $V(C)$, and $H_{1,E;r}(u)$
 is nonconstant on $V(C)$, $C$ being 
 the component  of $i'$ in $\cv\setminus\{i\}$,
  then there  exists $\witi Q\in \big(V_i\setminus\vo\big)\cap V(C)$  such that
  $H_{1,E;r}(u)(Q)=H_{1,E;r}(u)(\witi Q)$, and
   $\psi_i^{-1}(\witi Q)$  is an $E$-nonharmonic point  for $T_{i,E;r}(u)$.
 
\qua  iii) If  $E$ is a positive
  $r$-eigenform, $V_j$, $j=1,...,N$, has at least two multiple points, and 
 $\cv\setminus \{j\}$ is connected, then $\bar v_j$ has at
 least three  $E$-nonharmonic points. }
\smad

\qua i) follows from (2.4).

\qua ii) We use the notation of Lemma 3.5. If $\partial_1=\empty$,
as $Q\in V(C)\setminus \vo$, by Lemma 3.5 ii) $v$ is constant on $V(C)$, 
a contradiction. Thus $\partial_1\ne\empty$, and
$v(\witi Q)=v(Q)$ for some $\witi Q\in \partial_1$. By the definition
of $\partial_1$, $\witi Q\in V_h\cap V_{i'}$ for some 
$h\in C$ and $i'\notin C$, but as then $\{h,i'\}\in \cg$, by the definition of
$C$, we have $i'=i$. Thus,
$\witi Q$ satisfies the statement ii).
  
 \qua iii) If $\bar v_j$ attains its maximum at more than one point, the result is 
  trivial, as we can take the three extremum points. Suppose $\bar v_j$ attains
   its maximum at a unique point $P_{\bar j}$, so that
   
   $$\bar v_j(P_h)<\bar v(P_{\bar j}) \text{\  \  if\ } h\ne\bar j. \eqno (3.3)$$
   
    Then $P_j$ and $P_{\bar j}$ are
   two different $E$-nonharmonic points for $\bar v_j$. I will exhibit 
   a further $E$-nonharmonic point for $\bar v_j$ different from them.   
   Now, using notation of Lemma 3.5 with $C=\cv\setminus \{j\}$, 
    $u=\bar v_j$, and $v:=H_{1,E;r}(\bar v_j)$,
   we have

   $$T_{j,E;r}(\bar v_j)= l_j\bar v_j, \quad 0<l_j<1, \eqno (3.4)$$
   $$v\big(\psi_j(P_h)\big)=T_{j,E;r}(\bar v_j)(P_h)
   <v(P_h)\quad \forall\, h\ne j\, .\eqno (3.5)$$  
   Also, by (3.4) a point is $E$-harmonic for $\bar v_j$ if and only if
   it is $E$-harmonic for $T_{j,E;r}(\bar v_j)$. Moreover,
   
   $$  \partial_1=\psi_j(\partial_2) \eqno (3.6)
   $$
   
 In fact, 
 
 $$ \partial_1=\big\{\psi_j(P_h)\in V(\cv\setminus\{j\})
:  h\  E\text{-nonharmonic\  for\ }\bar v_j\big\}$$
$$=\big\{\psi_j(P_h):\#\big(st (\psi_j(P_h))\big)>1
:  h\  E\text{-nonharmonic\  for\ }\bar v_j\big\}$$
$$=\big\{\psi_j(P_h):h\ne j
:  h\  E\text{-nonharmonic\  for\ }\bar v_j\big\},$$
$$ \partial_2=\big\{P_h:h\ne j
:  h\  E\text{-nonharmonic\  for\ }\bar v_j\big\}.$$

   and (3.6) holds. By the previous formulas
   $P_{\bar j}\in \partial_2$ and $\psi_j(P_{\bar j})\in \partial_1\ne\empty$,
   and by (3.5) $v$ and does not attains
   its minimum on  $V(C)$ at points in $\partial_2$, hence by Lemma 3.5 i)
     $v$ attains its minimum on $V(C)$  at some point 
  $\psi_j(P_{\tilde j})\in\partial_1$.
   By   the hypothesis that $V_j$ has at least two multiple points,
   there exist $h\ne \bar j$ such that $\psi_j(P_h)\in V(C)$, thus,
   by (3.3) and (3.4)
    $$v\big(\psi_j(P_{\tilde j})\big)\le 
   v\big(\psi_j(P_h)\big)=l_j \bar v_j(P_h)
    <l_j \bar v_j(P_{\bar j})=  v\big(\psi_j(P_{\bar j})\big)$$
   
   hence $P_{\tilde j}\ne P_{\bar j}$,  and $P_j\ne P_{\tilde j}$
   as $\psi_j(P_{\tilde j})\in\partial_1$, and by formulas
   after (3.6), $P_{\tilde j}$ is $E$-nonharmonic point
   for $\bar v_j$. 
      \enpr

\bigskip
{\bf Remark 3.7.} We could expect that the converse of Lemma 3.6 i) holds
for nonconstant $u$, that is, if $Q\in V_i\setminus \vo$ is not $1$-point,
then $\psi_i^{-1}(Q)$ is an $E$-nonharmonic  point for $T_{i,E;r}(u)$.
 However, this is not in general true. For example, in the Vicsek set
 (Figure 4),
 if $u$ takes the value $0$ at $P_1$, $1$ at $P_3$, and ${1\over 2}$ at $P_2$
 and $P_4$,  and $r\equiv 1$, then $H_{1,E;r}(u)\equiv {1\over 2}$
on $V_2$ and $V_4$, as can be easily verified. Thus, for $i=5$,
$\psi_i^{-1}(Q)$ is an $E$-nonharmonic  point for $T_{i,E;r}(u)$, where
$Q$ is the point depicted in Figure 4, and the same for $Q'$. However, such a situation
appears to be rather unusual, that is, appears to hold only for
rather special fractals and/or functions $u$. \enpr

     \centerline{\bf 4. Generalizations}
     \medskip
     
     The argument in Theorem 3.1 can be used also for other kind of fractals. 
     In order to describe a  
     general situation
      I am now giving the following definition.
      {\it Put 
      
      $$S'_{i,B}:=\Big\{j\in\cu : st\big(\psi_i(P_j)\big)\not\subseteq B\cup \{ i\}\Big\},$$
      $$S_{i,B}:=S'_{i,B}\cup\Big (\{i\}\cap\{1,...,N\}\Big)$$ 
      
      when $i=1,...,k$, $B\subseteq \{1,...,k\}$.I say that $B\subseteq \{1,...,k\}$ is {\it \forcu} 
      by $A\subseteq \{1,...,N\}$ if $\# \big(S_{i,B}\cap A) < 2\Rightarrow i\in B$.}
      This definition is motivated by the following lemma.

      \medskip{
      \bf Lemma 4.1.} {\sl Suppose $r\in W$, $\alpha$ is  $+$ or $-$, 
      $u\in A^{\alpha} (E,E')$, $v=H_{1;E;r}(u)$, $E,E'$ are$\Lb_r$-eigenforms, 
      and $E$ has positive coefficients. Then the set $B:=\{i=1,...,k: v$ 
      is constant on $V_i\}$is \forcu\ by the set
       $A(\alpha):=\{j=1,...,N: u_{j;E;r}\in\ta^{\alpha} (E,E')\}$.} 
      
      \smallskip Proof. Assume for example $\alpha=+$ and put  $T_{i;E;r}=T_i$ . 
      If $P_j$ is an extremum point of $T_i(u)$ and  $j\notin S_{i,B}$,then $j\ne i$, thus
 $\psi_i(P_j)\in\von\setminus\vo$. Also,  $st\big(\psi_i(P_j)\big)\subseteq B\cup \{ i\}$, thus
       by Lemma 2.3 and the definition of $B$, $v$ is constant on $V_i$, 
       hence $T_i(u)$ is constant, 
       and $i\in B$. Thus, also using Prop. 2.1 iv and Lemma 2.2, we see that, if $i\notin B$,
then $S_{i,B}\cap A(+)$ contains at least two points, corresponding to the maximum and
         the minimum of $T_i(u)$. Therefore, if$\# \big(S_{i,B}\cap A(+)) < 2$, 
         then $i\in B$. \enpr
         
      \medskip
      As an example I now describe as the use of Lemma 4.1works for 
      {\it hexagonal pattern fractals}.Let $S=\{A_j: j=1,...,6\}$ be the set of the vertices
      of a regular hexagon of center $C$. We choose  the indices  fixing $A_1$ and putting 
      $A_{1+i}$ to be the vertex obtained  by  rotating $A_1$ of ${2i\pi\over 6}$
       around $C$.
      In particular, $A_j$ is opposite to $A_{j\pm 3}$.  Put $v_j=A_{j\pm 3}-A_j$, and
      $\wt=\big\{S+n v_j+m v_{j'}: n,m\in\bin\big\}$, where $0<|j-j'|<3$.
      Note that the set $\wt$ is independent of $j,j'$.
      
      \s, Based on $\wt$, we construct a fractal.
      We take $\cv=\{V_i:i=1,...,k\}\subseteq \widetilde \cv$, and$V=\{P_j: j=1,...,6\}$, 
      where every $P_j$ is a vertex  of some hexagon in $\cv$, and every $V_i$ amounts
       to $\psi_i(V)$ where $\psi_i$ has the form$\psi_i(x)=\alpha x+\beta_i$, $0<\alpha<1$. 
       We also require that the points $P_j$ are indexed so that $\psi_i(P_j)=A_j+B_i$
    for some $B_i\in \bre^2$. We assume that $\Psi=\{\psi_i: i=1,...,k\}$ generates a strongly 
     connected fractal $\cf$and that $V_i\subseteq co V$ for each $i=1,...,k$.In this case, 
     I will say that the fractal $\cf$ is a {\it hexagonal pattern} fractal or shortly h.p. fractal. 
     A usual example of this type is the Lindstr\o m snowflake. 
     Note that a h.p. fractal is never a tree. In  a h.p. fractal,
      I say that an$m$-plet $(i_1,...i_m)$
      with $i_l=1,...,k$ is a $j$-diagonal if it is a maximal set with the property that
      $\psi_{i_{l+1}}(P_j)=\psi_{i_l}(P_{j\pm 3})$.  I say that a $j$-diagonal is {\it long} 
      if $i_1=j, i_m\le N$, and that is short otherwise.

       \medskip
       {\bf Theorem 4.2.} {\sl On h.p. fractals we have at most one 
       normalized $\Lb_r$-eigenform, if $r\in W$.}
       
       \smallskip
        Proof. Suppose we have two  nonproportional $\Lb_r$-eigenforms $E$, $E'$. 
        Take $u$, $v$, $\alpha$ as in Lemma 4.1, and prove that $A(\alpha)$ has 
        at least four elements so that $\ta^+(E,E')\cap\ta^-(E,E')\ne\empty$, in contrast 
        to Prop. 2.1 iii. Note that $A(\alpha)$ has at least two elements. I will prove, 
        using Lemma 4.1, that if $A(\alpha)$ has two or three elements, we get a 
        contradiction.
        
        \s, Suppose first that $A(\alpha)$ has two elements $j$ and $j'$, 
        and take $u=u_j$. Then $v$ is constant on every $j$-diagonal $D$  
        not containing $j$. 
        Otherwise, let $D=(i_1,...,i_m)$ and let $\bar l$ be the minimum $l$ s
        o that $i_l\notin B$. 
        Then, $S_{i_{\bar l},B}\cap A(\alpha)\subseteq \{j'\}$, as
        $st\big(\psi_{i_{\bar l}}(P_j)\big)$ 
        amounts to$\{ i_{\bar l}, i_{\bar l-1 }\}$(or  to $\{ i_{\bar l }\}$,  if $\bar l=1$),
        thus is contained 
        in $B\cup\{i_{\bar l} \}$, and by hypothesis$i_l\ne j$ for all $l$, thus, by Lemma 4.1 
      $i_{\bar l}\in B$, a contradiction. 
      
      \s, Now, let $D=(i_1,...,i_m)$ be the $j$-diagonal containing 
      $j$ (hence $i_1=j$). As $v$ is constant on every $V_i$ with $i\notin D$, 
      by Lemma 2.3, 
      for every  $i\in D$, $v$ takes its extrema on $V_i$ at points either in $\vo$ 
or belonging to another cell in $D$. Note that on $V_j$,    $v$ is not constant and 
takes its minimum at 
  $\psi_j(P_j)$ by definition of $u_j$, hence we inductively get that on every $V_i$ 
  with $i\in D$,$v$ is nonconstant and takes its (strict) minimum  at $\psi_i(P_j)$ 
  and its (strict) 
      maximum  at $\psi_i(P_{j \pm 3})$. Also, $\psi_{i_m}(P_{j\pm 3})$ is a $1$-point, thus 
      by Lemma 2.3 again, and the fact that $v$ is nonconstant 
      on $V_{i_m}$, then$i_m=j\pm 3$.
       Therefore, $D$ is the long $j$-diagonal. By the strong connectedness and 
       the structure of h.p.  fractals, $V_{i_l}\cap V_i\ne\empty$ and
       $V_{i_{l+1}}\cap V_i\ne\empty$for some $1\le l<m$, $i\notin D$. 
       Let $Q\in V_{i_l}\cap V_i$,$Q'\in V_{i_{l+1}}\cap V_i$. On one hand, $v(Q)=v(Q')$ 
       as $v$ is constant on $V_i$, on the other, $v(Q')>v(\psi_{i_{l+1}}(P_j))=
       v(\psi_{i_{l}}(P_{j\pm 3}))>v(Q)$, a contradiction. 
       
       \s, Suppose instead  $A(\alpha)$ has 
       three elements $j$, $j'$ and $j''$, and take $u=u_j$ again in Lemma 4.1.
       We can assume, 
       up to a rotation, $j=1$, $j'=3$, and $j''\in \{2, 5, 6\}$.Note that $j\notin B$ as 
       $T_{j;E;r}(u_{j;E;r})$ is a positive multiple of $u_{j;E;r}$, thus is not constant. 
       Let $D$ be the $j'$ diagonal containing $j$, $D=(i_1,...,i_m)$, 
       and we have $i_m=j$. Also, 
       let $\bar l$ be the minimum $l$ such that $i_l\notin B$.Note that 
       $j,j'\notin S'_{i_{\bar l},B}$ 
       as $st\big(\psi_{i_{\bar l}}(P_j)\big)=\{i_{\bar l}\}$,
       $st\big(\psi_{i_{\bar l}}(P_{j'})\big)=\cases \{i_{\bar l}\}&\ \text{if}\ \bar l=1\\
      \{ i_{\bar l},i_{\bar l-1}\}&\ \text{if}\ \bar l>1\endcases$.
       Now, if $\bar l\ne m$, then $j,j'\notin \{i_{\bar l}\}$,  
       therefore$S_{i_{\bar l},B}\cap A(\alpha)\subseteq\{j''\}$, and by Lemma 4.1
        $i_{\bar l}\in B$, 
       a contradiction. If, instead, $\bar l= m$, then $i_{\bar l}=j$, and $j,j'\notin S'_{j,B}$. 
       Moreover, also $j''\notin S'_{j,B}$, as if $j''$ is $2$ or $6$, 
       then $st\big(\psi_{j}(P_{j''})\big)=\{j\}$, and if $j''=5$, then
       the considerations for $j'$ are also valid for $j''$.
       Hence, $S_{j,B}\cap A(\alpha)\subseteq\{j\}$, and by Lemma 4.1 again, 
       $j\in B$, a contradiction. \enpr

       \medskip I now describe  examples of fractals in which we get the uniqueness of 
       the eigenform, implicitly using in fact not Lemma 4.1 but a variant of it, which I will 
       not explicitly state.

          \bigskip
          {\bf Remark 4.4.} I call {\it tree-diagonal} a long diagonal $D$ 
          in the sense of Remark 3.2 
          such that there exists no loop containing an element of $D$. 
          In the previous examples 
           the nonuniqueness occurs only on quasi-tree fractals(cf. Remark 3.2). 
           However, a slightly deeper analysis of such examples as well as of the arguments
           used in the proof suggests that the uniqueness of the normalized 
           eigenform could be 
           related to the fact that there exists at most one  tree-diagonal rather than that
           the fractal is quasi-tree.  \enpr
           
           \bigskip

\end